\pdfoutput=1
\newif\ifpreprint
\preprinttrue

\ifpreprint
  \documentclass[12pt, 
                 numbers=noenddot,
                 headings=normal, 
                 DIV16, BCOR10mm,
                 enabledeprecatedfontcommands]{scrartcl}
  \usepackage{natbib}
  \usepackage{graphicx}
  \usepackage{color}
  \usepackage{tikz}
\else
  \documentclass[preprint,12pt,authoryear]{elsarticle}
  \journal{Advances in Water Resources}
\fi

\usepackage{amsmath,amssymb}
\usepackage{stmaryrd}
\usepackage{subfig}
\usepackage{units}
\usepackage{booktabs}
\usepackage{xspace}
\usepackage[hyphens]{url}
\usepackage[linktocpage,breaklinks]{hyperref}
\usepackage{lineno}

\newcommand{\absperm}{\mathbb{K}}       
\newcommand{\flux}{{q}}             
\newcommand{\normal}{\mathbf{n}}        
\newcommand{\pressure}{{p}}         
\newcommand{\source}{{q}}           

\newcommand{\dirichlet}{\text{\footnotesize D}} 
\newcommand{\neumann}{\text{\footnotesize N}}       

\newcommand{\K}{\absperm}

\newcommand{\neumannflux}{\flux_\neumann}

\newcommand{\pdir}{\pressure_\dirichlet}

\newcommand{\globalDomain}{\mathcal{D}}
\newcommand{\pDn}{{\partial\globalDomain_{\neumann}}}
\newcommand{\pDd}{{\partial\globalDomain_{\dirichlet}}}
\newcommand{\vecu}{\mathbf{u}}

\DeclareMathOperator{\grad}{\mathbf{grad}}

\DeclareMathOperator{\Div}{div}
\renewcommand{\div}{\Div}

\newcommand{\transp}{{\top}}

\newcommand{\Dumux}{DuMu$^\text{x}$\xspace}

\newcommand{\m}{\mathrm{m}}
\newcommand{\f}{\mathrm{f}}

\newcommand{\nsub}{\mathrm{n}}
\newcommand{\tsub}{\mathrm{t}}

\newcommand{\base}{p}

\newcommand{\scalarPermeability}{k}
\newcommand{\tensorPermeability}{\boldsymbol{K}}
\newcommand{\Km}{{\tensorPermeability_{\m}}}

\newcommand{\kft}{\scalarPermeability_{\f,\tsub}\,}
\newcommand{\kfn}{\scalarPermeability_{\f,\nsub}\,}
\newcommand{\aperture}{\varepsilon}

\newcommand{\ie}{i.e. }

\newcommand{\eg}{e.g. }

\newcommand{\unitTensor}{{\mathbf{I}}}

\newcommand{\mean}[2]{\left\{\!\!\left\{ {#1} \right\}\!\!\right\}_{#2}}

\newcommand{\jump}[2]{\left\llbracket {#1} \right\rrbracket_{#2}}

\graphicspath{{./pics/}}

\newcounter{BenchmarkCounter}
\stepcounter{BenchmarkCounter}

\ifpreprint

\else
\begin{document}

\begin{frontmatter}
\fi

\title{Benchmarks for single-phase flow in fractured porous media}

\ifpreprint
\newcommand*{\affaddr}[1]{{\normalsize #1}}
\newcommand*{\affmark}[1][*]{\textsuperscript{#1}}
\newcommand*{\email}[1]{\texttt{#1}}
\author{Bernd Flemisch\affmark[1], Inga Berre\affmark[2], Wietse Boon\affmark[2], Alessio Fumagalli\affmark[2],\\Nicolas Schwenck\affmark[1], Anna Scotti\affmark[3], Ivar Stefansson\affmark[2],\\and Alexandru Tatomir\affmark[4]\\
{\normalsize\affmark[1]Department of Hydromechanics and Modelling of Hydrosystems, University of Stuttgart,}\\
{\normalsize Pfaffenwaldring 61, 70569 Stuttgart, Germany, \texttt{bernd@iws.uni-stuttgart.de}}\\
{\normalsize\affmark[2]Department of Mathematics, University of Bergen, All\'{e}gaten 41, 5007 Bergen, Norway}\\
{\normalsize\affmark[3]Laboratory for Modeling and Scientific Computing MOX, Politecnico di Milano,}\\
{\normalsize p.za Leonardo da Vinci 32, 20133 Milano, Italy}\\
{\normalsize\affmark[4]Department of Applied Geology, Geosciences Center, University of G\"ottingen,}\\
{\normalsize Goldschmidtstrasse 3, 37077 G\"ottingen, Germany}}
\else
\author[1]{Bernd Flemisch}
\author[2]{Inga Berre}
\author[2]{Wietse Boon}
\author[2]{Alessio Fumagalli}
\author[1]{Nicolas Schwenck}
\author[3]{Anna Scotti}
\author[2]{Ivar Stefansson}
\author[4]{Alexandru Tatomir}
\address[1]{Department of Hydromechanics and Modelling of Hydrosystems, University of Stuttgart, Pfaffenwaldring 61, 70569 Stuttgart, Germany}
\address[2]{Department of Mathematics, University of Bergen, All\'{e}gaten 41, 5007 Bergen, Norway}
\address[3]{Laboratory for Modeling and Scientific Computing MOX, Politecnico di Milano, p.za Leonardo da Vinci 32, 20133 Milano, Italy}
\address[4]{Department of Applied Geology, Geosciences Center, University of G\"ottingen, Goldschmidtstrasse 3, 37077 G\"ottingen, Germany}
\fi

\ifpreprint

\begin{document}
\maketitle
\fi

\begin{abstract}
This paper presents several test cases intended to be benchmarks for
numerical schemes for single-phase fluid flow in fractured
porous media.  A number of solution strategies are compared, including a
vertex and a cell-centered finite volume method, a non-conforming embedded
discrete fracture model, a primal and a dual extended finite element
formulation, and a mortar discrete fracture model.  The proposed benchmarks test
the schemes by increasing the difficulties in terms of network geometry,
\eg intersecting fractures, and physical parameters, \eg low and high
fracture-matrix permeability ratio as well as heterogeneous fracture
permeabilities.  For each problem, the results presented by the participants are
the number of unknowns, the approximation errors in the porous matrix and in the
fractures with respect to a reference solution, and the sparsity and condition
number of the discretized linear system. All data and meshes used in this study
are publicly available for further comparisons.
\end{abstract}

\ifpreprint
\else
\begin{keyword}
fractured porous media \sep discretization methods \sep benchmark
\end{keyword}

\end{frontmatter}
\fi


\section{Introduction}
\label{sec:introduction}
In porous-media flow applications, the domains of interest often contain geometrically anisotropic inclusions and strongly discontinuous material coefficients that can span several orders of magnitude. If the size of these heterogeneities is
small in normal direction compared to the tangential directions,
these features are called fractures. Fractures can act both as conduits and barriers and affect flow patterns severely. Target applications concerning fractured porous-media systems in earth sciences include groundwater resource management, renewable energy storage,
recovery of petroleum resources, radioactive waste reposition, coal bed methane migration in mines, and geothermal energy production.

The analysis and prediction of flow in fractured porous media systems are
important for all the aforementioned applications. Many different conceptual and
numerical models of flow in fractured porous-media systems can be found in the
literature. Even though fractured porous-media systems have been of interest to
modelers for a long time, they still present challenges for simulators.
During the last 70 years, different modeling approaches have been developed and
gradually improved. Comprehensive reviews can be found in
\citet{berkowitz2002characterizing,dietrich2005flow,hoteit2008numerical,neumann2005review,sahimi2011flow,singhal2010applied}.
Roughly, the fractured porous media systems are classified in two broad
categories: discrete fracture-matrix  (DFM) models and continuum fracture
models. Within this paper, we will only consider DFM models.

The DFM models consider flow occurring in both the fracture network and the
surrounding rock matrix. They account explicitly for the effects of individual
fractures on the fluid flow. An efficient way to represent fractures in DFMs
is the hybrid-dimensional approach, \eg
\citet{helmig_multiphase_1997,firoozabadi_control-volume_2004,karimi2004efficient,martin2005modelling,reichenberger_mixed-dimensional_2006},
where fractures in the geometrical domain are discretized with elements of co-dimension one with respect
to the dimension of the surrounding matrix, such as one-dimensional elements in two-dimensional
settings. The aforementioned classical DFM approaches all rely on matching fracture and
matrix grids in the sense that a fracture element coincides geometrically
with co-dimension-one mesh entities, \ie faces of matrix grid elements.
In addition to the classical models, several so-called non-conforming DFM models have been developed
in recent years, such as EDFM \citep{moinfar2013development,hajibeygi2011hierarchical}, XFEM-based approaches
\citep{dangelo2012mixed,Schwenck:2015:DRF,huang2011use}, or mortar-type methods \citep{frih2012modeling}.

Benchmarking represents a methodology for verifying, testing and comparing the
modeling tools.  Various codes have been developed by academic institutions or
companies based on different conceptual, mathematical, and numerical models.
Even though benchmarking studies are increasing in all fields
of engineering and workshops have been organized around specific problems (\eg
\cite{class_benchmark_2009}), there are still only a limited number of
studies of this type in the field of geoscience.
Some are related to a specific application and are flexible as to how
the problem is modeled in terms of assumptions regarding the physics and the selection of
the domain, see
\cite{dahle_model-oriented_2010,nordbotten_uncertainties_2012,caers_special_2013,kolditz_thermo-hydro-mechanical-chemical_2015}.
Others \citep{deDreuzy2013, caers_special_2013}, like ours, focus on the
comparison of numerical schemes.
One of the common requirements when selecting the test problems for
comparing numerical schemes is that they allow the examination of the
capabilities of each of the compared methods. Therefore, our benchmark study
proposes a set of problems starting from simple geometries and then gradually
increasing the geometrical complexity.  The test problems are specifically
selected to make clear distinctions between the different methods.

The main focus of this work is to use existing and new
computational benchmarks for fluid flow in fractured porous media to
compare several DFM-based numerical schemes in a systematic way.
We would also like to invite the scientific community to follow up on this
study and evaluate further methods by means of the proposed benchmarks.
In order to facilitate this, the paper is accompanied by grid and result
files in the form of a Git repository at
\url{https://git.iws.uni-stuttgart.de/benchmarks/fracture-flow}.

The remainder of this paper is organized as follows.
In Section \ref{sec:model}, we formulate the model problem in terms of
the partial differential equation to be solved.
The participating DFM models are described in Section \ref{sec:methods}.
The central Section \ref{sec:benchmarks} proposes the benchmarks and
compares the results of the different methods.
Finally, Section \ref{sec:summary} concludes with a summary and outlook.

\section{The model problem}
\label{sec:model}
We are considering an incompressible single-phase flow through a porous medium,
assumed to be described by Darcy's law, resulting in the governing sytem of equations
\begin{subequations}\label{eq:strong}
\begin{align}
 \vecu &= - \K \grad\pressure, \label{eq:darcy} \\
\div \vecu &= \source, \label{eq:massbalance}
\end{align}
in an open bounded domain $\mathcal{D} \subset \mathbb{R}^N$, subject to boundary conditions
\begin{alignat}{3}
 \pressure &= \pdir & & \text{ on } \pDd, \label{eq:bcdir} \\
 \vecu\cdot\normal &= \neumannflux && \text{ on } \pDn, \label{eq:bcneu}
\end{alignat}
\end{subequations}
with $\partial\mathcal{D} = \overline{\pDd \cup \pDn}$ and $\pDd \cap \pDn =
\emptyset$. In equations (\ref{eq:strong}), $\mathbf{u}$ denotes the macroscopic
fluid velocity whereas $\K$ and $\pressure$ stand for absolute permeability and
pressure.

Let us assume that $\globalDomain$ contains several fractures, that all together
constitute a single domain $\Gamma$ of spatial dimension $N$ such that $\Gamma
\subset \globalDomain$, which is a possibly unconnected, open subset of
$\globalDomain$. The surrounding porous rock, namely, the remaining part of
$\globalDomain$, is called $\Omega = \globalDomain \setminus \overline{\Gamma}$.
Assuming that the fracture aperture $\aperture$ at each point of $\Gamma$ is small compared to
other characteristic dimensions of the fractures, the full-dimensional domain
$\Gamma$ can be reduced to the $(N{-}1)$-dimensional fracture network $\gamma$.
This reduction involves modeling choices resulting in different hybrid-dimensional
problem formulations that form the basis for the methods presented in the following section.

\section{Participating discretization methods}
\label{sec:methods}
Within this section, the discretization methods participating
in this benchmark study are described. The purpose of this article is the comparison of
well-known, established and/or at least published methods.
Therefore, only the
most significant aspects of each method are summarized. We do not show a
comparison against analytical solutions here. The analysis of the methods and theoretical results such as proofs of optimal convergence can be found in the corresponding references. A summary of all participating methods is provided in Table \ref{tbl:methods}. In the sequel, we will denote with \textit{d.o.f.} the
degrees of freedom associated to a specific method. We indicate also the type of
conformity required to the computational grid with respect to the fractures and the
assumption that the pressure is considered continuous across the fractures.
\renewcommand{\arraystretch}{1.1}
\begin{table}[hbt]
\centering
\begin{tabular}{|l|c|c|c|c|}\hline
\textbf{method} & \textbf{d.o.f.}
& \textbf{frac-dim} & \textbf{conforming} & $p$\textbf{-cont.} \\\hline
Box-DFM & $p$ (vert) & dim-1 & yes & yes \\\hline
CC-DFM & $p$ (elem) & dim-1 & yes & no \\\hline
EDFM & $p$ (elem) & dim-1 & no & yes \\\hline
mortar-DFM & $p$ (elem), $\boldsymbol{u}$ (faces) & dim-1 & geometrically & no \\\hline
P-XFEM & $p$ (vert) & dim-1 & no & no \\\hline
D-XFEM & $p$ (elem), $\boldsymbol{u}$ (faces) & dim-1 & no & no \\\hline
MFD & $p$ (faces) & dim & geometrically & no \\\hline
\end{tabular}
\caption{Participating discretization methods.
}
\label{tbl:methods}
\end{table}
\renewcommand{\arraystretch}{1.0}

\subsection{Vertex-centered, continuous-pressure, conforming lower-dimensional DFM (Box-DFM)}
\label{sec:boxDFM}
The lower-dimensional representation of fractures allows easier mesh generation
in comparison to the equi-dimensional approach, as it circumvents the appearance
of very small elements when discretizing the interior of the fracture (\ie,
within the fracture width).
The conforming mesh generation algorithm honors the geometrical
characteristics of the fracture system. Conform meshing implies that the
fractures are discretized with a set of line elements (in a 2D domain) that
are also the edges of the triangular finite elements.

\begin{figure}[hbt]
\centering
  \includegraphics[width=1\linewidth]{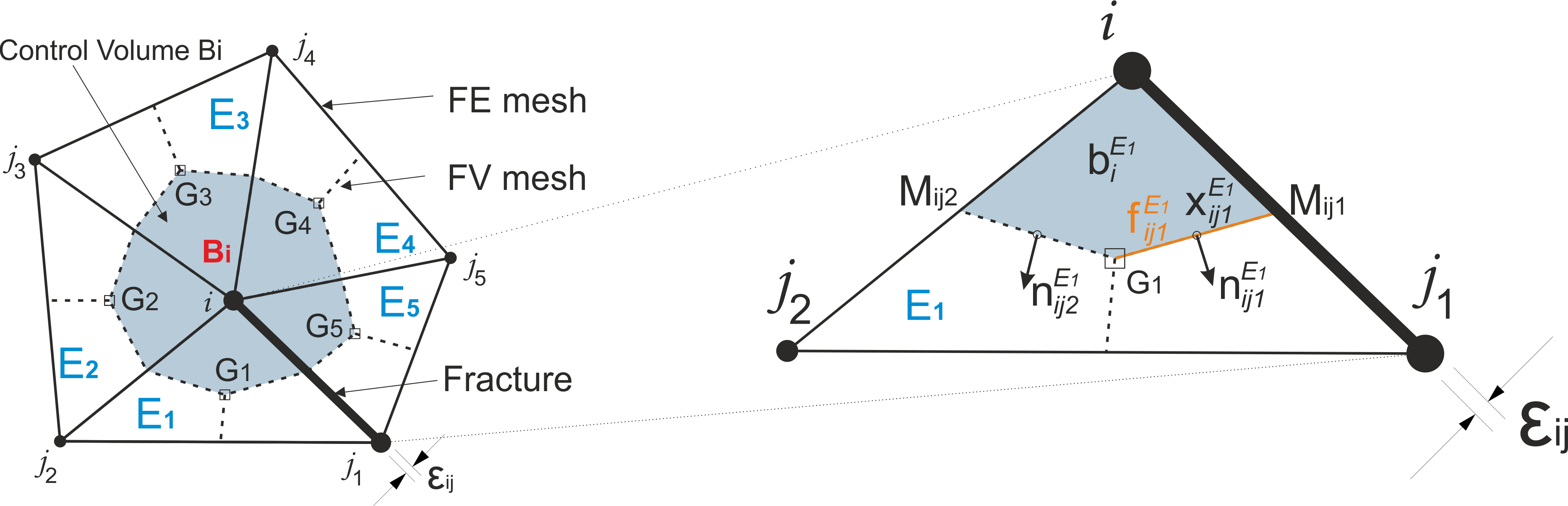}
\caption{Conceptual representation of the Box-DFM method:  (left-hand side) The dual  finite element and finite volume mesh from which the control volume \textit{B$_i$} around node \textit{i} is created. Node $i$ is surrounded by nodes \{$j_1$, $j_2$, $j_3$, $j_4$, $j_5$\}, where segment $ij_1$ represents both a fracture and a shared FE edge; (right-hand side) Sub-control volume (SCV) $b_i^{E_1}$  in element $E_1$ has barycenter $G_1$ and the mid-points of the edges $ij_1$ and $ij_2$ are $M_{ij1}$, respectively $M_{ij2}$. The SCV face $f_{ij1}^{E1}$ is the segment $\overline{G_1 M_{ij1}}$ which contains the integration point $x_{ij1}^{E1}$ where the normal vector $\textbf{n}_{ij1}^{E1}$ is applied.}
\label{fig:BoxDFMConcept}
\end{figure}

The spatial discretization in Box-DFM is performed with the Box method, a vertex-centered finite-volume method proposed in, \eg
\citet{helmig_multiphase_1997} which combines the advantages of finite element
and finite volume grids, allowing unstructured grids and guaranteeing {a}
locally conservative {scheme} (\cite{ reichenberger_mixed-dimensional_2006}).
Figure \ref{fig:BoxDFMConcept} illustrates a two-dimensional representation of
the dual-grid with  two  finite elements $E_1$ and $E_5$  sharing the same edge
($ij_1$) that represents a lower-dimensional fracture with the aperture
$\aperture_{ij}$.
The main characteristic in terms of the fractured system is that the pressure is required to be continuous,
in particular in those vertices whose control volumes overlap both fracture and matrix regions.

The Box-DFM method used for this paper is implemented in the open-source numerical simulator \Dumux.
A detailed description of the conceptual, mathematical and numerical model and code implementation is published in \cite{Tatomir12}.
The Box-DFM simulation code used for the benchmark studies is publicly available under
\url{https://git.iws.uni-stuttgart.de/dumux-pub/Flemisch2016a.git}.

\subsection{Cell-centered, discontinuous-pressure, conforming DFM (CC-DFM)}
The control volume finite difference method uses a two-point flux approximation (TPFA) based on the cell-center pressure values for the evaluation of the face fluxes,
and is a widely applied and standard method for simulation of flow in porous media. The domain is partitioned with fractures coinciding with the interior faces between matrix cells just as described in Section \ref{sec:boxDFM}.
The flux over the face between matrix cells $i$ and $j$ is approximated by
\begin{equation}\label{tpfa}
\vecu_{ij} = T_{ij}(\pressure_{i}-\pressure_{j}),
\end{equation}
where $\pressure_{i}$ and $\pressure_{j}$ are the pressures in the neighboring cells and $T_{ij}$ is the face transmissibility, computed as the harmonic average of the two half transmissibilities corresponding to the face and the two cells. The half transmissibility of cell-face pair $i$ is in turn given as
\begin{equation} 
\alpha_{i}  =  \dfrac{A_{i}  \normal_{i}^\transp \K_{i} }{\mathbf{d}_{i}^\transp
\mathbf{d}_{i} } \cdot \mathbf{d}_{i},
\end{equation}
where $A_{i}$ and $\normal_i$ are the area and unit normal vector of the face,
$\K_i$ is the permeability assigned to the cell and $\mathbf{d}_{i}$ is the distance vector from cell center to face centroid.

\begin{figure}[hbt]
\centering
\subfloat[]{
  \includegraphics[width=0.35\linewidth]{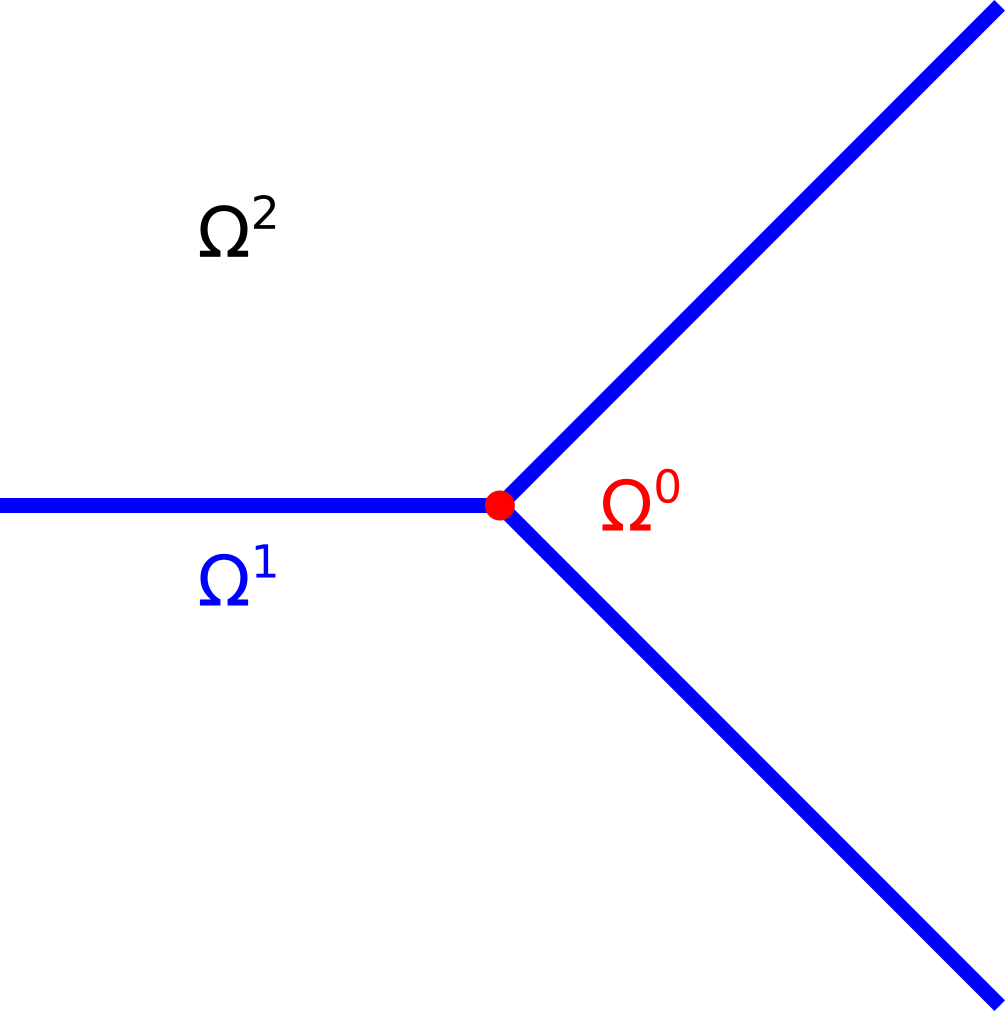}
}
\subfloat[]{
  \includegraphics[width=0.35\linewidth]{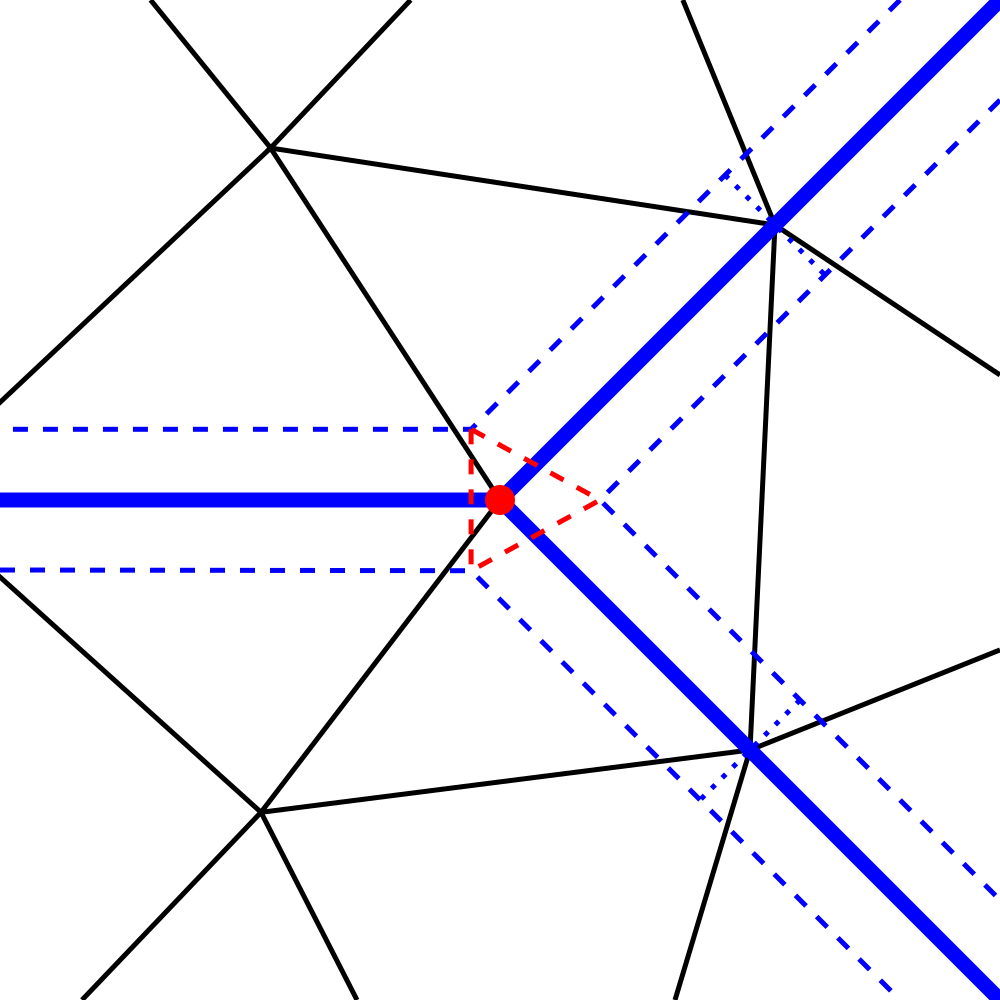}

}
\caption{(a) Conceptual decomposition of the domain according to element dimension with the matrix depicted in black, fractures in blue and their intersections in red. (b) The computational domain of the CC-DFM. Dashed lines are faces of the fracture cells.}
\label{fig:Domain_dimensions}
\end{figure}

In addition to the unknowns given at the centroids of the matrix cells, unknowns are associated to the centroids of the fracture cells. The fracture cells are associated with apertures, which multiplied with the length give the volume of these cells. The aperture is also used to construct hybrid faces for the matrix-fracture interfaces. These faces, parallel to the fracture but displaced half an aperture to either side, enable us to compute the half transmissibilities between the fracture cell and the matrix cells on the two sides. These faces are indicated by the dashed blue lines in  Figure \ref{fig:Domain_dimensions}, where the computational domain is superimposed on the geometrical grid.
The result is a hybrid grid with fractures which are lower dimensional in the grid, but equidimensional in the computational domain at the cost of a small matrix volume error corresponding to the overlap of the matrix cells with the fracture cells.

Following the method proposed by \cite{karimi2004efficient}, the intermediate fracture intersection cell drawn with dashed red lines in Figure \ref{fig:Domain_dimensions} is removed, leading to direct coupling of the fracture cells neighbor to the intersection. The purpose of this is both to obtain a smaller condition number and to avoid severe time-step restrictions associated with small cells in transport simulations. To each new face between cell \emph{i} and \emph{j}, face transmissibilities are assigned, calculated using the star delta transformation as described in \cite{karimi2004efficient}:
\begin{equation}\label{halftransmissibility}
T_{ij}  =  \dfrac{\alpha_{i}  \alpha_{j} }{\sum\limits_{k=1}^n{\alpha_{k}}},
\end{equation}
with $n$ denoting the number of fracture cells meeting at the intersection. As this elimination disregards all information on the permeability of the intersection, it should be used with caution in cases of crossing fractures of different permeability. We encounter this feature in \ref{sec:anna}, and include results both with and without the elimination for one of the test cases presented in that section.

Inspired by the CC-DFM method by \cite{karimi2004efficient} presented above, a method based on the multi-point flux approximation has also been developed \cite{sandve2012efficient}. The MPFA variant of the method reduces errors associated with the TPFA approach for grids that are not close to K-orthogonal, and avoids errors related to the splitting of the fluxes in the star-delta transformation. We refer to \cite{sandve2012efficient} for a thorough comparison of the TPFA and MPFA CC-DFM approaches. The implementation of both methods is available in the open-source Matlab Reservoir Simulation Toolbox \url{http://www.sintef.no/projectweb/mrst/}.

\subsection{Continuous-pressure, non-conforming embedded DFM (EDFM)}
Recently, non-conforming methods for the treatment of lower-dimensional
fractures have been developed, for example in
\cite{moinfar2013development,moinfar2011comparison,hajibeygi2011hierarchical},
to avoid the time-consuming construction of complex matrix grids which
explicitly represent the fractures.  They are mostly used in the context of
single and multi-phase flow simulations for petroleum engineering applications
and require the normal fracture permeability to be orders of magnitude higher
than the matrix permeability, as in the case of enhanced reservoir exploitation and
fractures stimulation.
In this field of applications corner-point grids are normally employed to describe the geological layers,
\eg different rock type, of the reservoir. An adaptation of {such} computational
grids to the fractures could be unaffordable for real cases.
The numerical method belongs to the family of two-point schemes, where a
one-to-one connection between the degrees
of freedom is considered  through the transmissibility concept (\cite{Eymard2000}). References on the embedded discrete
fracture method (EDFM) can be found, for example, in
\cite{Li2008,Panfili2013,Moinfar2014,Panfili2014,AraujoCavalcanteFilho2015,Fumagalli2015b}.

\begin{figure}[hbt]
    \centering
    \includegraphics[width=0.45\linewidth]{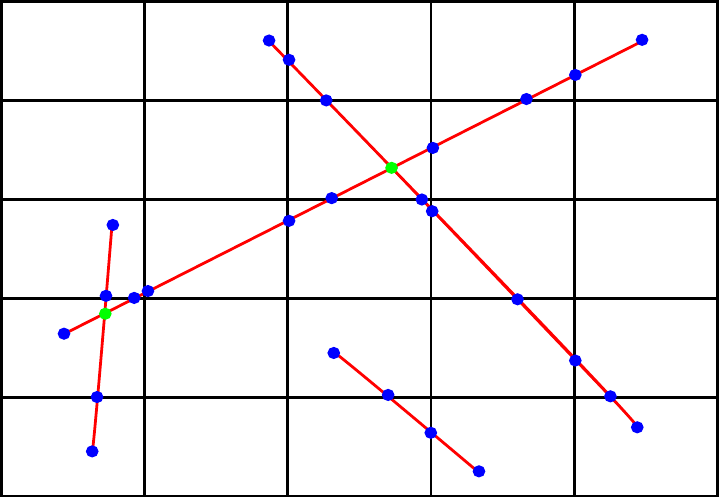}
    \caption{Example of meshes, for both fractures and rock matrix, suited for
    EDFM. The rock matrix is considered as a background mesh. Each fracture
    cell is represented by two blue dots and the green dots are the
    non-matching intersection among fractures.}%
    \label{fig:edfm_grid_intro}
\end{figure}

In practice, the meshes of the fractures are generated on top of the rock grid
so that each rock cell cut by fractures contains exactly one fracture cell per
fracture.  Intersections between fractures are computed without affecting the
creation of the grids of fractures and rock and used to compute
approximate transmissibilities between different fracture cells. See Figure
\ref{fig:edfm_grid_intro} as an example.
A degree of freedom that represents a pressure or a saturation value is assigned to
each matrix cell and to each fracture cell.  This means that transmissibilities
between matrix and fracture cells, as well as those between different fracture
cells, need to be computed. We compute the transmissibility between a fracture cell and a
matrix cell $T_{fm}$ and the half-transmissibility $T_i$ between two intersecting
fracture cells (related to the fracture $i$) through the following approximate
expressions
\begin{gather*}
    T_{fm} = A \dfrac{ \normal_f^\transp \K \cdot \normal_f}{d_{f,m}}
    \quad \text{and} \quad
    T_{i} = s \dfrac{k_i \aperture_i}{d_{i,s}}.
\end{gather*}
Here $A$ is the measure of the fracture cell in the current rock cell,
$\normal_f$ is the normal of the fracture cell and $d_{f,m}$ is an average
distance between the fracture cell and the matrix cell, see \cite{Li2008}.  For the
fracture-fracture transmissibility, $s$ indicates the measure of the intersecting
segment, $k_i$ the scalar permeability of the fracture, $\aperture_i$ the aperture and
$d_{i,s}$ is the average distance between the fracture cell and the intersecting
segment.  The standard harmonic average is considered to compute the
transmissibility between the two fracture cells.  Standard formulae for
fracture-fracture as well as matrix-matrix transmissibilities are computed by means of a
two-point flux approximation. It is worth to notice that the recent extension
 of EDFM called Projection-based EDFM (pEDFM), proposed in \cite{Tene2016},  is
 also able to handle low permeable fractures. Finally, even if the proposed
 benchmark cases
are two-dimensional the method can be extended to three dimensions without any
additional constraints.

\subsection{Cell-centered, discontinuous-pressure, geometrically-conforming mortar DFM (mortar-DFM)}
The key concept behind the mortar-DFM, as described more thoroughly in \citet{Boon2015Robust}, is the idea that fractures can be considered as interfaces between different sub-domains. This has been explored previously by \citet{martin2005modelling,frih2012modeling}, among others. In this context, it is interesting to consider domain decomposition techniques such as the mortar method to model flow through the fractured porous medium.

\begin{figure}[hbt]
    \centering
    \includegraphics[width=0.35\linewidth]{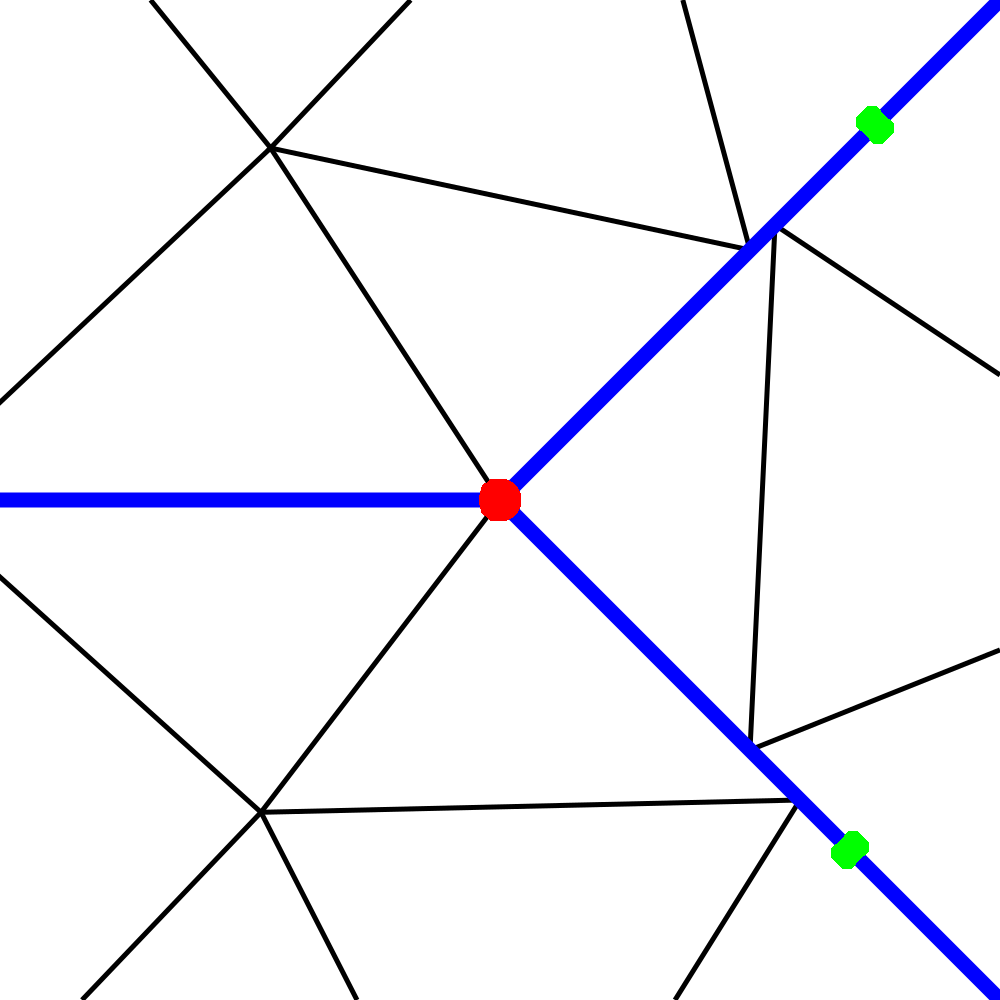}
    \caption{The mortar-DFM allows for non-matching grids along fracture interfaces. Fracture and matrix flows are coupled using a mortar variable, defined on a coarser grid (green dots).}%
    \label{fig:mortar_grid_intro}
\end{figure}

The mortar method is generally used to couple equations in different sub-domains by introducing a so-called mortar variable, defined on the interface. In case of modeling fracture flow, a well-explored choice of the mortar variable is the fracture pressure \citep{martin2005modelling}. The method considered here, however, uses as the mortar variable the flux, which leads to a stronger sense of mass conservation for flows between the matrix and fractures. One of the main advantages of the close relationship to mortar methods is the capability to handle non-matching grids. In particular, two sub-domains bordering a fracture can be meshed independently on both sides,
as illustrated in Figure \ref{fig:mortar_grid_intro}.
The difficulty in mesh generation is then relieved significantly since only the geometry of the fractures needs to be respected.

By construction, the mortar-DFM is applicable to problems in arbitrary dimensions. The governing equations in the matrix and the fractures (as well as fracture intersections in 3D) are identical and thus all fractures, intersections and tips are handled in a unified manner. Consequently, although only two-dimensional problems are considered in this case study, the discretization scheme is not at all limited to the presented benchmark problems and 3D cases can easily be considered.

With the use of mixed finite elements, mass is conserved locally in the matrix,
fractures, and fracture intersections. The flux $\vecu$ and pressure $\pressure$
are modeled using the lowest-order Raviart-Thomas elements and piecewise
constants, \ie $\mathbb{RT}_0-\mathbb{P}_0$. Additionally, the mortar variable is given by piecewise constants on a separately generated, lower-dimensional, mortar grid. This grid matches with the surrounding grids in case of matching grids and is coarser otherwise \citep{Boon2015Robust}. The resulting mixed finite element formulation is a saddle-point problem, which may be challenging to solve numerically. To relieve this, the flux variables may be eliminated through hybridization, which leads to a less expensive scheme containing solely the cell-center pressures.

Two implementations of the method have been developed, both of which are used in this benchmark study. The first version, implemented in MATLAB, has the capability of handling non-matching grids along fractures for two-dimensional problems. It is well-suited for simpler geometries, containing relatively few fractures, such as those considered in Benchmarks 1-3. The second version has been implemented for 3D problems and higher-order spaces on matching grids using the open-source finite element library FEniCS \citep{LoggMardalEtAl2012a}. This version is more efficient for complex cases such as Benchmark 4.

\subsection{Discontinuous-pressure, non-conforming primal XFEM (P-XFEM)} The
primal XFEM method participating in this benchmarking study is described in detail
in \cite{schwenck2015thesis}, see also
\cite{Flemisch:2016:RXA,Schwenck:2015:DRF}.  The method is based on the
hybrid-dimensional problem formulation investigated in
\cite{martin2005modelling}, where conditions for the coupling between fracture
and matrix are derived:
\begin{subequations}\label{eq:interface}
\begin{align}
\mean{\vecu_\m\cdot \normal}{\gamma} &= \kfn/\aperture \jump{\pressure_\m}{\gamma} \\
\xi_0 \jump{\vecu_\m \cdot \normal}{\gamma}
&= \kfn/ \aperture \left(\mean{\pressure_\m}{\gamma} - \pressure_\f \right)
\end{align}
\end{subequations}
Here, the subscripts ``m'' and ``f'' indicate matrix and fracture quantities,
while $\mean{\cdot}{\gamma}$ and $\jump{\cdot}{\gamma}$ denote the average and
the jump of a matrix quantity over the fracture $\gamma$, respectively.

The coupling conditions \eqref{eq:interface} can be used to define a source term
for the fracture flow problem, while they yield an interface problem for the
matrix domain. For the discretization of this interface problem, the
methodology presented in \cite{hansbo2002unfitted} is used, which amounts to
applying the eXtended Finite Element Method (XFEM).  Together with an independent
standard discretization of the lower-dimensional fracture problem, this yields a
hybrid-dimensional, non-conforming primal XFEM-based method.
The XFEM space is built enriching the standard Lagrangian $\mathbb{P}_1$
(or $\mathbb{Q}_1$ for quads) finite-element spaces, whose degrees of freedom are located at the vertices of the full-dimensional
grid of the matrix $\Omega$ and the lower-dimensional grid of the fracture
$\gamma$.
\begin{figure}[hbt]
    \centering
    \includegraphics[width=0.45\linewidth]{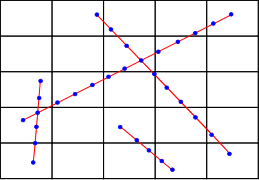}
    \caption{Example of meshes, for both fractures and rock matrix, suited for
    P-XFEM. The fracture
    grid vertices are indicated by the blue dots.}%
    \label{fig:xfem_grid_intro}
\end{figure}
A representative example of matrix and fracture grids is illustrated in Figure
\ref{fig:xfem_grid_intro}. Unlike the EDFM method, see Figure \ref{fig:edfm_grid_intro}, the fracture grid vertices can be placed arbitrarily
without taking into account the matrix grid. On the other hand, the method requires matching fracture
branch grids in the form of vertices placed at the fracture intersections. In particular,
special care has to be taken of intersecting and immersed fractures
\citep{Schwenck:2015:DRF}.

The method is implemented on top of the DUNE framework
\cite{bastian2008generic2} and the discretization module DUNE-PDELab
\cite{bastian2010generic}. For the enrichment of the finite-element spaces in
the context of XFEM, the modules DUNE-Multidomain and DUNE-Multidomaingrid are
employed \cite{Muething:2015:FFM}. The simulation code for the XFEM approach and
for the benchmarks studied here is publicly available under
\url{https://git.iws.uni-stuttgart.de/dumux-pub/Flemisch2016a.git}.

\subsection{Discontinuous-pressure, non-conforming dual XFEM (D-XFEM)}

\begin{figure}[hbt]
\centering
  \includegraphics[width=0.8\linewidth]{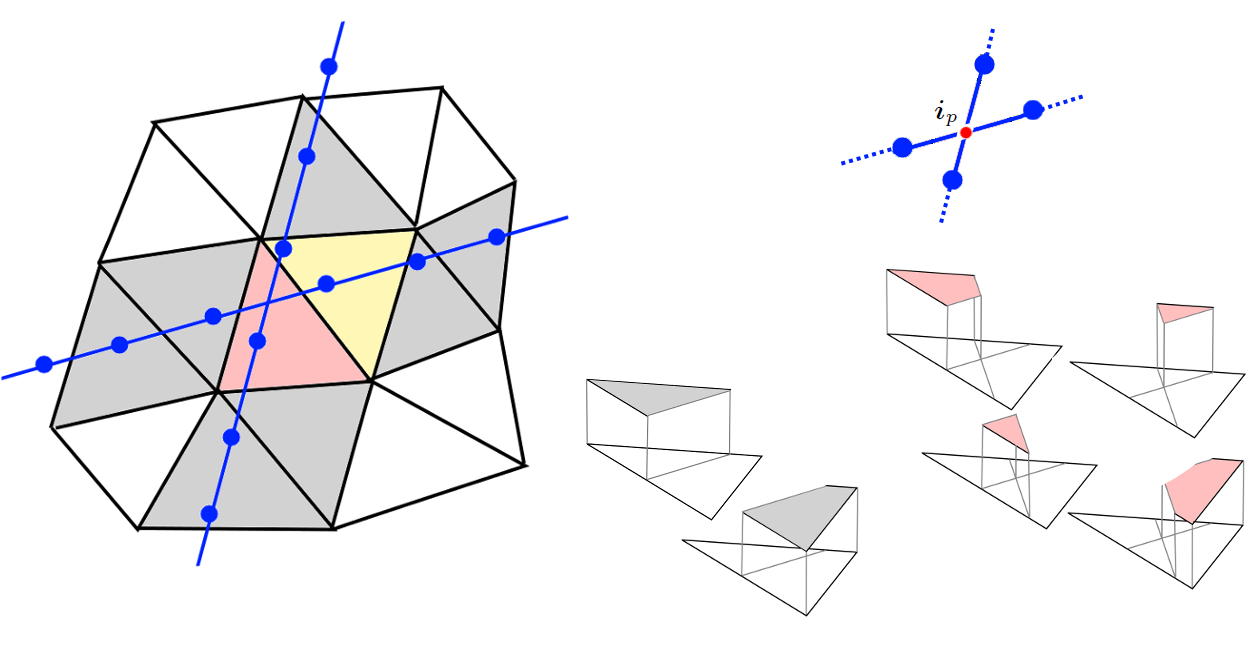}
\caption{A portion of the grid cut by two fractures: in the two dimensional case
they can split the elements in two (grey), three (yellow), or four (red)
independent parts, where the restrictions of the basis functions are defined.
The fracture grids are irrespective of the bulk grid and of each other, \ie the intersection point $i_p$ is not a point of the grid.}
\label{fig:Domain_DXFEM}
\end{figure}

The dual XFEM method participating in his benchmark is based on
\cite{dangelo2012mixed}. The method, originally derived for a domain cut by one
fracture, was further developed in \cite{formaggia2012crossings},
\cite{fumagalli2013ogst} to account for intersecting fractures with different
permeabilities.  The same equations and coupling conditions as for the primal
XFEM are used, but in a dual formulation where Darcy law and mass conservation
give rise to a saddle-point problem for the fluid mean velocity and pressure,
both in the fracture and in the surrounding medium. Moreover, unlike the previous method, this method employs triangular/tetrahedral grids. The usual lowest order
$\mathbb{RT}_0-\mathbb{P}_0$ pair for velocity and pressure is
enriched following  \cite{hansbo2002unfitted} in the elements of the porous
medium cut by a fracture, or in the elements of a fracture at the intersection
with other fractures. Indeed, triangular/tetrahedral grids are arbitrarily cut
by triangulated lines/surfaces in 2D and 3D respectively. These surfaces can, in turn,
intersect each other in a non-conforming way, as shown in Figure \ref{fig:Domain_DXFEM}.

In the current implementation of the method no special enrichment is added in
the bulk elements containing the fracture tips. Instead, fractures are
artificially extended up to the boundary of the domain, and in the extension we prescribe the same permeability of the surrounding porous medium to obtain a ``virtual'' fracture with no effects on the flow.

The method has been implemented on the basis of the Getfem++ library,
\url{http://download.gna.org/getfem/html/homepage/},
which provides support for the computation of the intersections and the
quadrature on sub-elements thanks to an interface with QHull,
\url{http://www.qhull.org/}.

\subsection{Reference Solutions calculated with mimetic finite differences (MFD)}
The reference solutions are computed on very fine grids that discretize both matrix
and fractures by full-dimensional triangular or quadrilateral elements.
A mimetic finite difference method, see \cite{BrezziLipnikovSimonici:2005:FMF,Flemisch:2008:NIM},
is used to discretize problem \eqref{eq:strong}. The
method is employed as it is implemented in \Dumux 2.7 \cite{flemisch2011dumux}.
In particular, a mixed-hybrid approach is used to transform the discrete saddle point problem
in terms of cell pressures and face fluxes into a symmetric positive definite formulation
with face-pressure degrees of freedom.

\section{Benchmark Problems}
\label{sec:benchmarks}
This is the main section which compares the methods described above by means of four
benchmark cases.
First, in Section \ref{sec:hydrocoin}, we present a well established benchmark
for groundwater flow from \cite{hydrocoin} that contains two crossing,
highly permeable fractures and a non-straight top surface. The second benchmark case, considered in Section \ref{sec:geiger_example}, is based on \cite{geiger2011novel} and shows a regular fracture network.
After that, a small but complex fracture network exhibiting immersed fractures and intersections at different angles
is investigated in Section \ref{sec:anna}.
Finally, a case synthesized from a real application is considered in Section \ref{sec:realistic}.

For each benchmark case, a description of the computational domain is provided, including
boundary conditions,
the geometrical information about the corresponding fracture network
and the associated material parameters such as aperture and permeability.
For some of the cases, the reference solution on the complete domain is visualized.
This is followed by illustrations of the grids used by the participating methods.
Since the methods pose different requirements there, the grid could be chosen
arbitrarily for each method, provided that the number of grid cells or vertices
is roughly the same.
If a reference solution is available (Benchmarks 1--3), the results of the different methods
are compared by evaluating the errors with respect to the reference in the matrix domain
as well as in the fracture network, indicated by $err_\text{m}$ and $err_\text{f}$, respectively.
The errors are calculated according to the formulas
\begin{align*}
err_\text{m}^2 &= \frac 1{|\Omega|(\Delta p_\text{ref})^2} \sum_{f = K_\text{ref} \cap K_\text{m}}
|f| \left(p_\text{m}|_{K_\text{m}} - p_\text{ref}|_{K_\text{ref}}\right)^2, \\
err_\text{f}^2 &= \frac 1{|\gamma|(\Delta p_\text{ref})^2} \sum_{e = K_\text{ref} \cap K_\text{f}}
|e| \left(p_\text{f}|_{K_\text{f}} - p_\text{ref}|_{K_\text{ref}}\right)^2,
\end{align*}
where $|\Omega|$ and $|\gamma|$ indicate the size of the matrix and fracture domain, respectively,
and $\Delta p_\text{ref} = \max_\globalDomain p_\text{ref} - \min_\globalDomain p_\text{ref}$.
The sum is taken over all intersections of (full-dimensional) elements $K_\text{ref}$
of the grid employed for the reference solution with full-dimensional
matrix elements $K_\text{m}$ in case of $err_\text{m}$ and lower-dimensional fracture elements
$K_\text{f}$ in case of $err_\text{f}$.
Moreover, the densities and condition numbers of the resulting linear system matrices are provided.
In addition to that, a comparison is performed by means of plots along specific lines through
the domain for some benchmark cases. Each case is concluded by a short discussion of the results.

\clearpage
\subsection{Benchmark \theBenchmarkCounter: Hydrocoin}
\label{sec:hydrocoin}

\begin{figure}[hbt]
\centering
\def\svgwidth{0.85\linewidth}
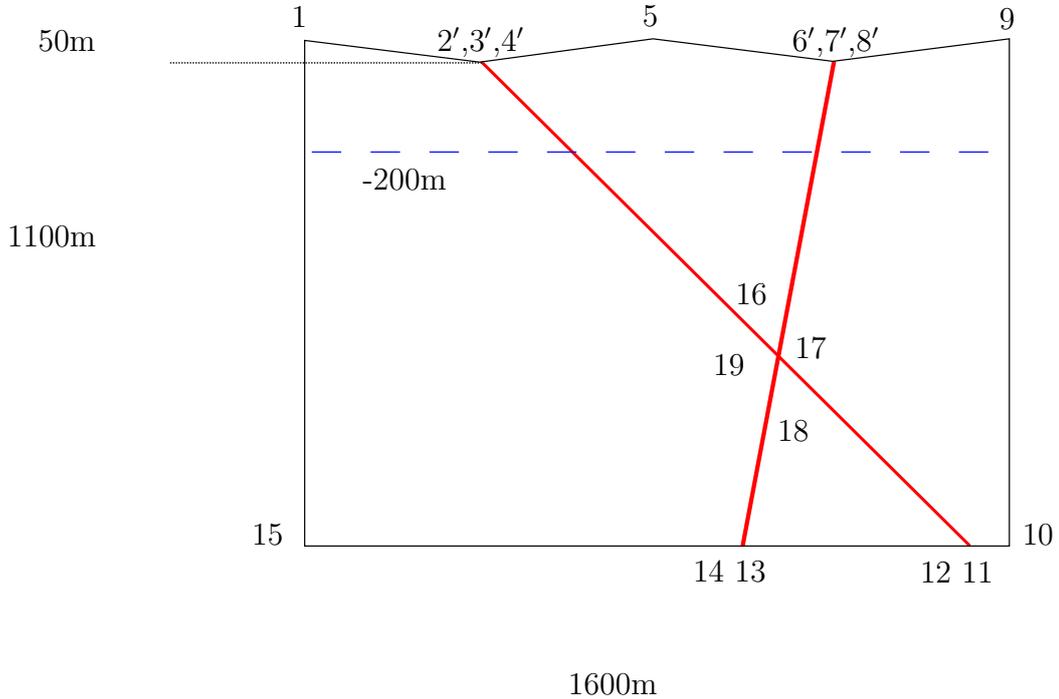
\caption{Geometry of the modeled domain of the Hydrocoin test case 2,
\cite{hydrocoin}. Modified node locations are indicated by numbers superscripted
with $\mbox{}^\prime$. Boundary conditions are hydraulic head on top and Neumann
no-flow on the other three sides of the domain.}
\label{fig:hydrocoin_setup}
\end{figure}

Within the international {Hydrocoin} project, (\cite{hydrocoin}), a benchmark for
heterogeneous groundwater flow problems was presented. The domain setup is shown
in Figure \ref{fig:hydrocoin_setup}.
We point out that we have slightly modified the original domain such that
equi-dimensional and hybrid-dimensional models can be run on exactly the same domain. This allows for an easier comparison of the solution
values over the whole domain.
The exact modifications are described in \ref{sec:hydrocoin_mod}.

For this case, we keep the original formulation in terms of the piezometric head 
and the hydraulic conductivity instead of pressure and permeability.
In particular,
the boundary conditions are Dirichlet~piezometric head on the top boundary and
Neumann~no flow on the other three boundaries. The hydraulic conductivity is $\unitfrac[10^{-6}]{m}{s}$ in the
fracture zones  and $\unitfrac[10^{-8}]{m}{s}$ in the rock matrix respectively.

\begin{figure}[hbt]
\centering
\subfloat[Box-DFM]{
  \includegraphics[width=0.315\textwidth]{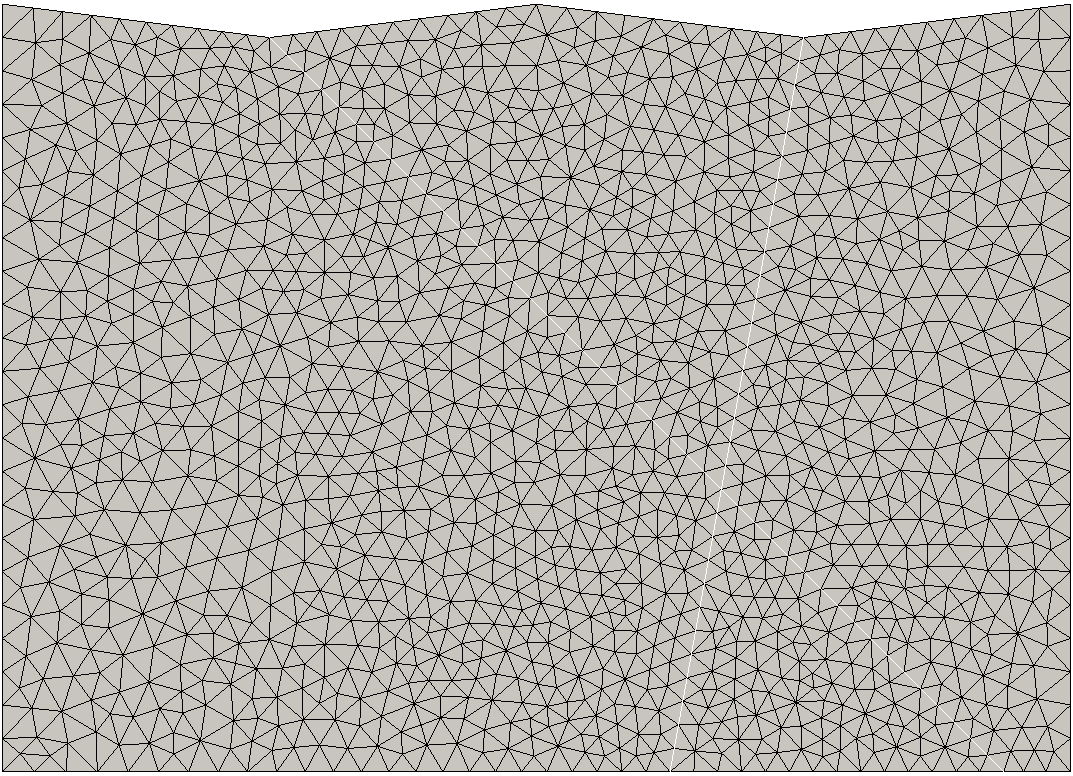}
}
\subfloat[CC-DFM]{
  \includegraphics[width=0.315\textwidth]{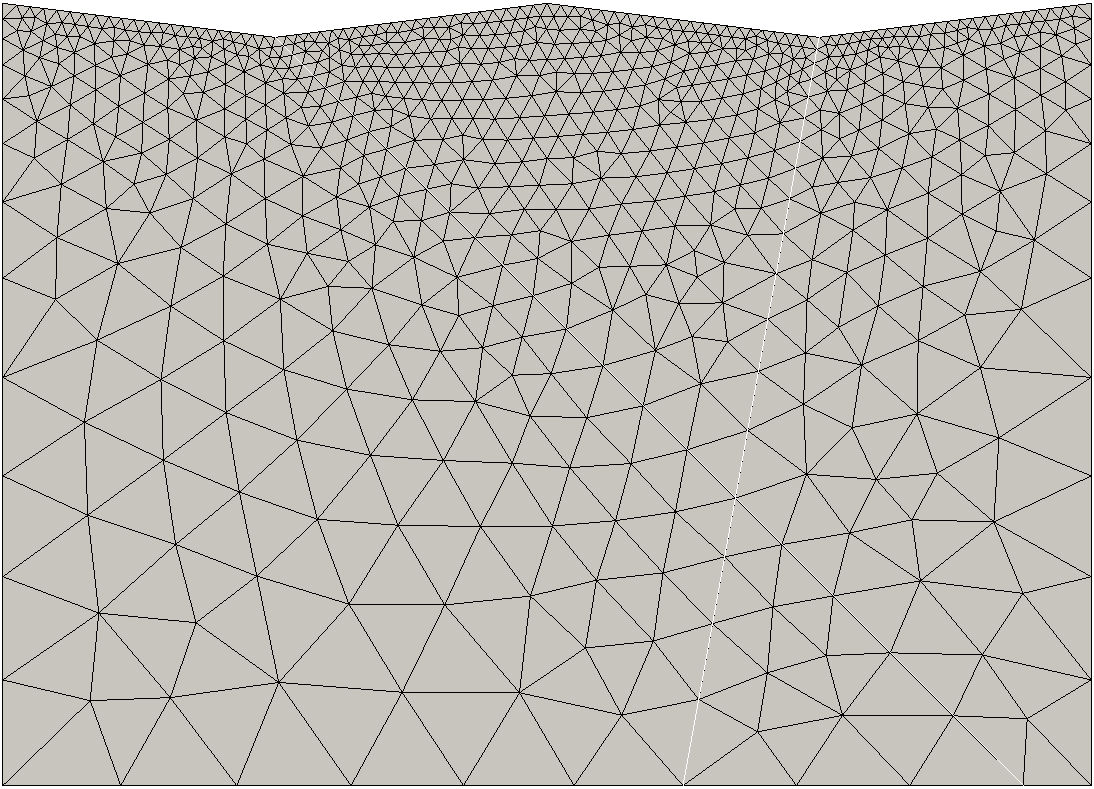}
}
\subfloat[EDFM]{
  \includegraphics[width=0.315\textwidth]{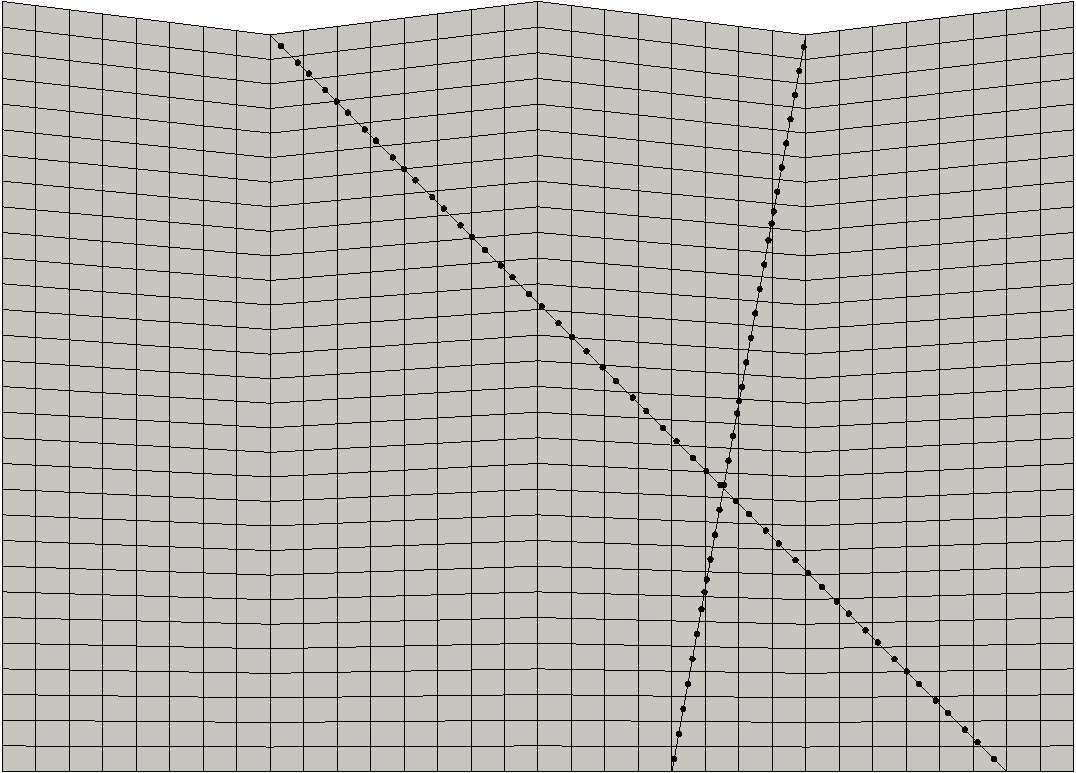}
}
\\
\subfloat[mortar-DFM]{
  \includegraphics[width=0.315\textwidth]{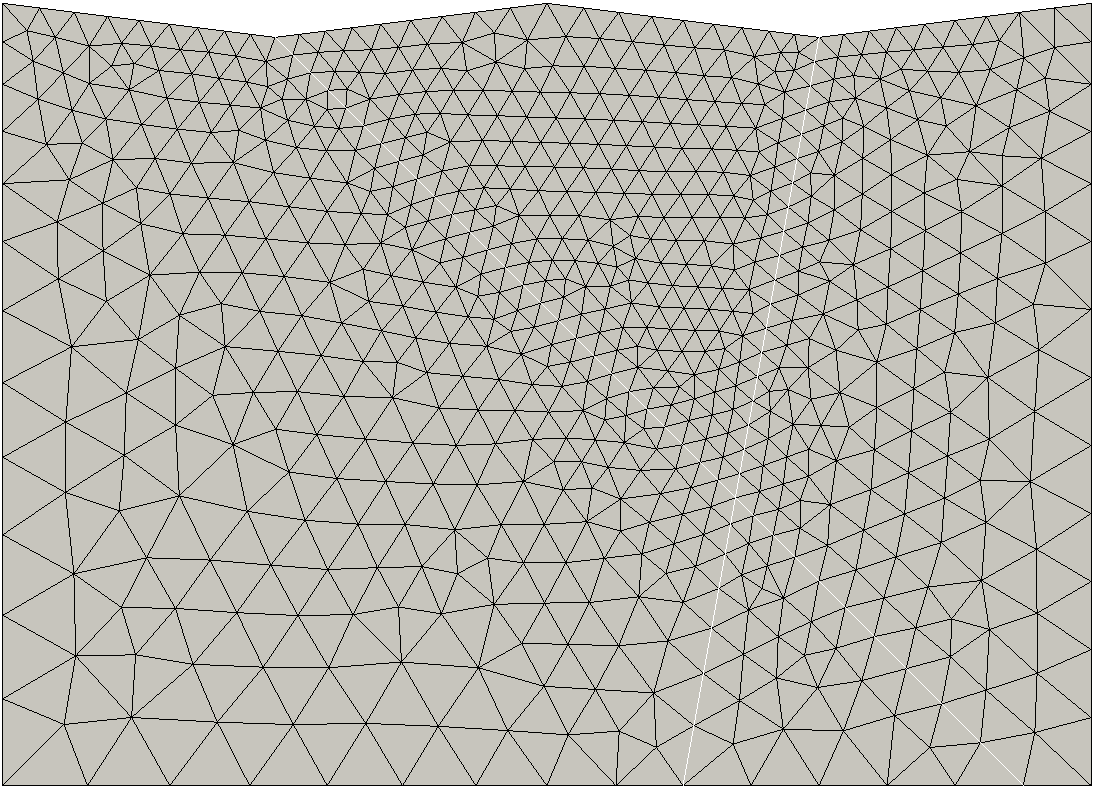}
}
\subfloat[P-XFEM]{
  \includegraphics[width=0.315\textwidth]{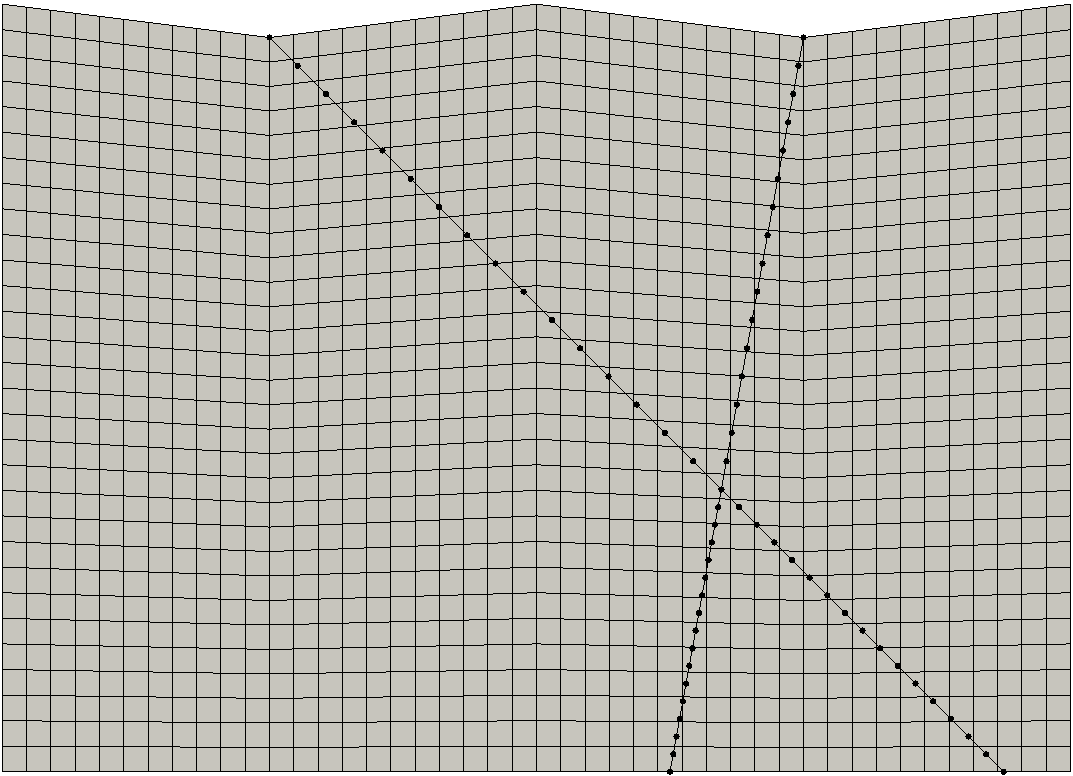}
}
\subfloat[D-XFEM]{
  \includegraphics[width=0.315\textwidth]{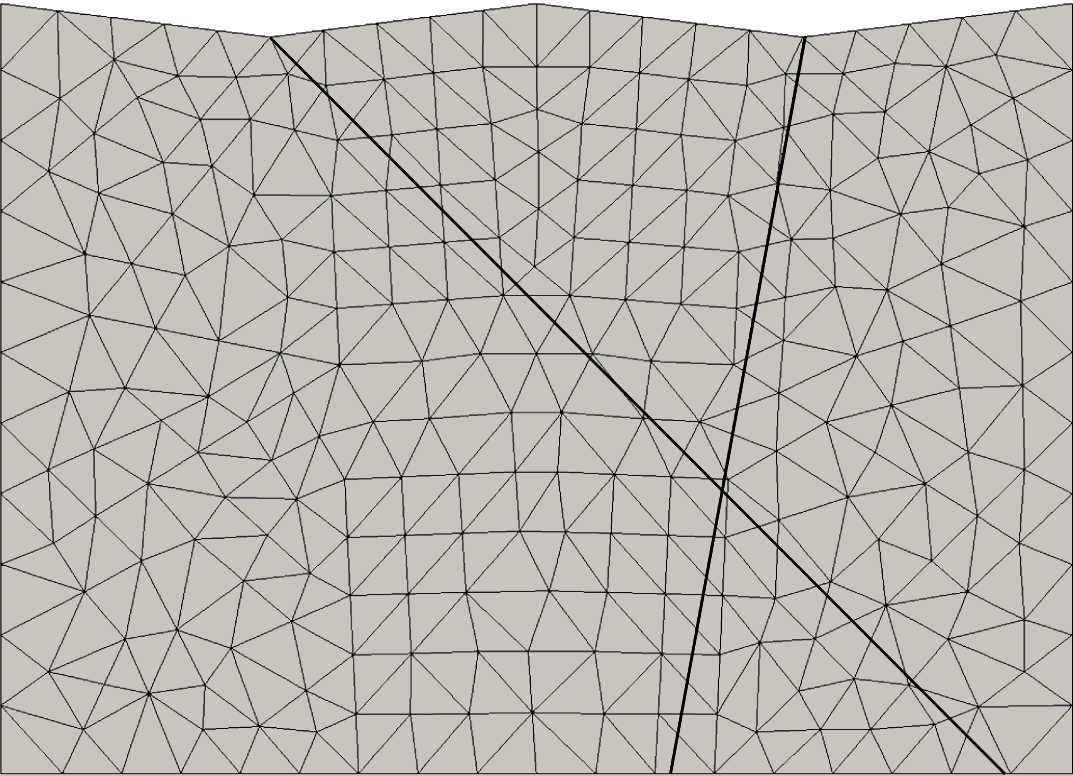}
}
\caption{Benchmark \theBenchmarkCounter: the grids used by the different methods.}
\label{fig:hydrocoin_grids}
\end{figure}

Table \ref{tbl:hydrocoin_grids} lists the number of degrees of freedom, matrix
elements and fracture elements for all the participating methods.
\renewcommand{\arraystretch}{1.1}
\begin{table}[hbt]
\centering
\begin{tabular}{|l|c|c|c|}\hline
\textbf{method} & \textbf{d.o.f.}
& \textbf{matrix elements} & \textbf{fracture elements} \\\hline
Box-DFM    & 1496 & 2863 triangles & 74 \\\hline
CC-DFM     & 1459 & 1416 triangles & 43 \\\hline
EDFM       & 1044 & 960 quads & 84 \\\hline
mortar-DFM & 3647 & 1384 triangles & 63 \\\hline
P-XFEM     & 1667 & 1320 quads & 68 \\\hline
D-XFEM     & 3514 & 1132 triangles & 160 \\\hline
MFD        & 889233 & 424921 mixed & 19287 \\\hline
\end{tabular}
\caption{Grids for Benchmark \theBenchmarkCounter.}
\label{tbl:hydrocoin_grids}
\end{table}
\renewcommand{\arraystretch}{1.0}
The corresponding grids are visualized in Figure \ref{fig:hydrocoin_grids}.

\begin{figure}[hbt]
\centering
\includegraphics[width=0.9\textwidth]{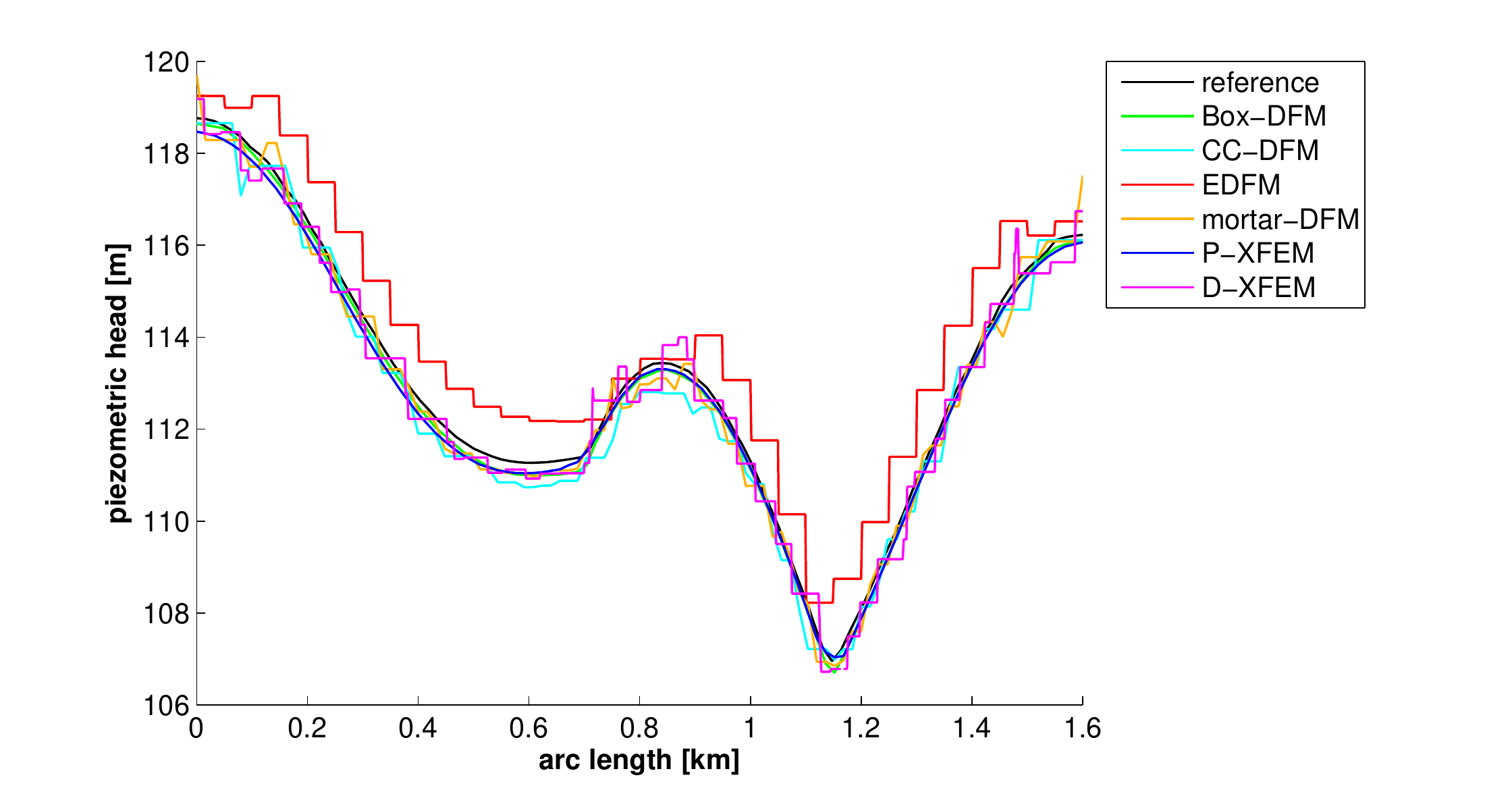}
\caption{Benchmark \theBenchmarkCounter: head values along a horizontal line at a depth of $\unit[200]{m}$.}
\label{fig:hydrocoin_200m}
\end{figure}

The original benchmark shows the piezometric head distribution along five horizontal lines
through the modeled domain. Here, we first show in Figure \ref{fig:hydrocoin_200m} the plot at a depth of $\unit[200]{m}$, as indicated by the dashed line in Figure \ref{fig:hydrocoin_setup}.
All participating methods show a good agreement with the reference solution. Only the EDFM method is
a bit off. We remark that the plots for the methods employing cell-wise constant solution values
exhibit staircase-like patterns corresponding to these values.

Table \ref{tbl:hydrocoin_error} lists the discretization errors for the different methods, particularly, the error for the matrix domain and the one along the two fractures.
\renewcommand{\arraystretch}{1.1}
\begin{table}[hbt]
\centering
\begin{tabular}{|l|c|c|c|c|}\hline
\textbf{method} & $err_\text{m}$ & $err_\text{f}$
& \textbf{nnz/size$^2$} & $\|\cdot\|_2$\textbf{-cond} \\\hline
Box-DFM    & 9.3e-3 & 3.3e-3 & 4.5e-3 & 5.4e3 \\\hline
CC-DFM     & 1.1e-2 & 1.1e-2 & 2.7e-3 & 3.5e4\\\hline
EDFM       & 1.5e-2 &  8.3e-3 & 4.7e-3 & 3.9e4 \\\hline
mortar-DFM & 1.0e-2 & 7.2e-3 & 1.5e-3 & 9.0e12\\\hline
P-XFEM     & 1.2e-2 & 3.2e-3 & 6.5e-3 & 2.7e9 \\\hline
D-XFEM     & 1.2e-2 & 6.9e-3  & 1.7e-3 & 6.2e12\\\hline
\end{tabular}
\caption{Discretization errors and matrix characteristics for Benchmark \theBenchmarkCounter.}
\label{tbl:hydrocoin_error}
\end{table}
\renewcommand{\arraystretch}{1.0}
Moreover, it provides the density of the associated matrix and its condition number for each method.
The uniform behavior exhibited in Figure \ref{fig:hydrocoin_200m} is reflected by the error values.
Especially the errors in the matrix domain are within very narrow bounds, while the fracture errors
show a larger variation. The densities of the matrices are also close together.
Remarkably high differences can be observed for the matrix condition numbers. While the ones for
Box-DFM, CC-DFM and EDFM are on the order of $10^4$, the one for P-XFEM is five orders and the ones
for mortar-DFM and D-XFEM are even seven orders of magnitude larger, due to their saddle-point nature.

\clearpage
\stepcounter{BenchmarkCounter}
\subsection{Benchmark \theBenchmarkCounter: Regular Fracture Network}
\label{sec:geiger_example}
\begin{figure}[hbt]
    \centering
    \def\svgwidth{0.7\linewidth}
    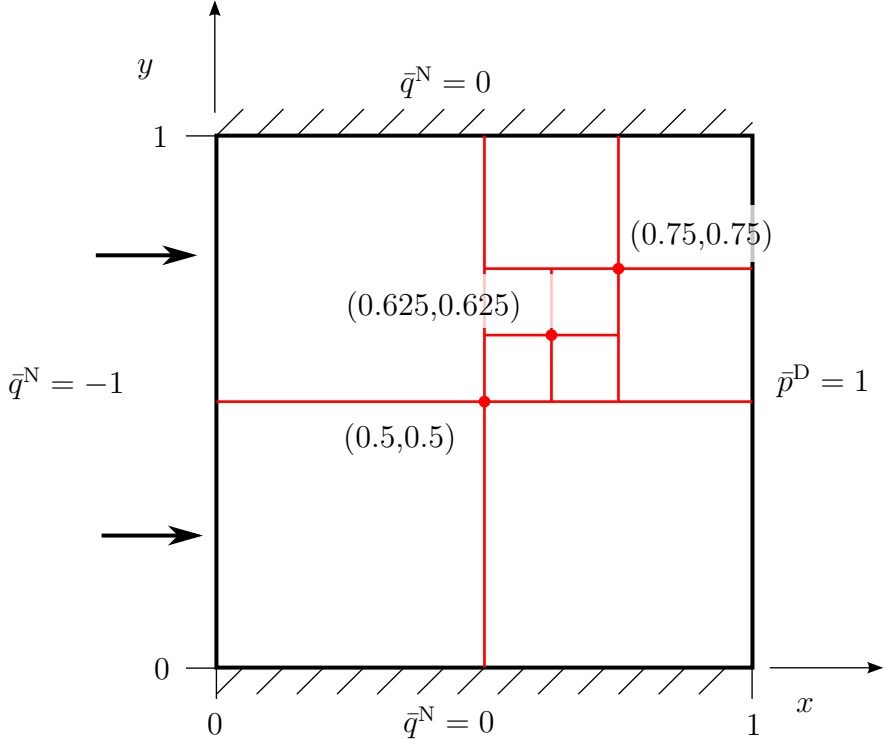
    \caption{Benchmark \theBenchmarkCounter: Domain and boundary conditions.}
\label{fig:geiger_domain}
\end{figure}

This test case is based on an article presenting a new dual continuum model, \cite{geiger2011novel}, with slightly modified boundary conditions and material properties. The computational domain including the fracture network and boundary conditions is shown in Figure \ref{fig:geiger_domain}.
The matrix permeability is set to $\Km=\unitTensor$, all fractures have a uniform aperture $\aperture=10^{-4}$.
For the fracture permeability we consider two cases: a highly conductive network with $\kfn=\kft=10^4$, as worked out in Subsection \ref{sec:geiger_example1}, and a case with blocking fractures by setting $\kfn=\kft=10^{-4}$,
as described in Subsection \ref{sec:geiger_example2}.
The reference solutions are computed on a grid which resolves every fracture with 10 elements in its normal direction and becomes coarser away from the fractures. It has a total of 1175056 elements.

\begin{figure}[hbt]
\centering
\subfloat[Box-DFM]{
  \includegraphics[width=0.315\textwidth]{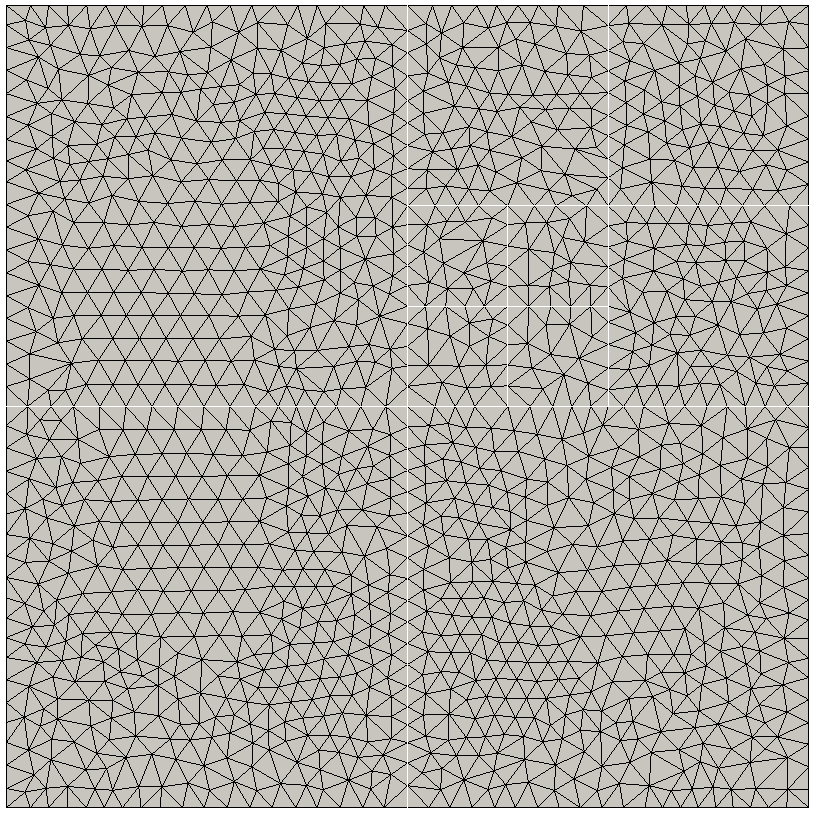}
}
\subfloat[CC-DFM]{
  \includegraphics[width=0.315\textwidth]{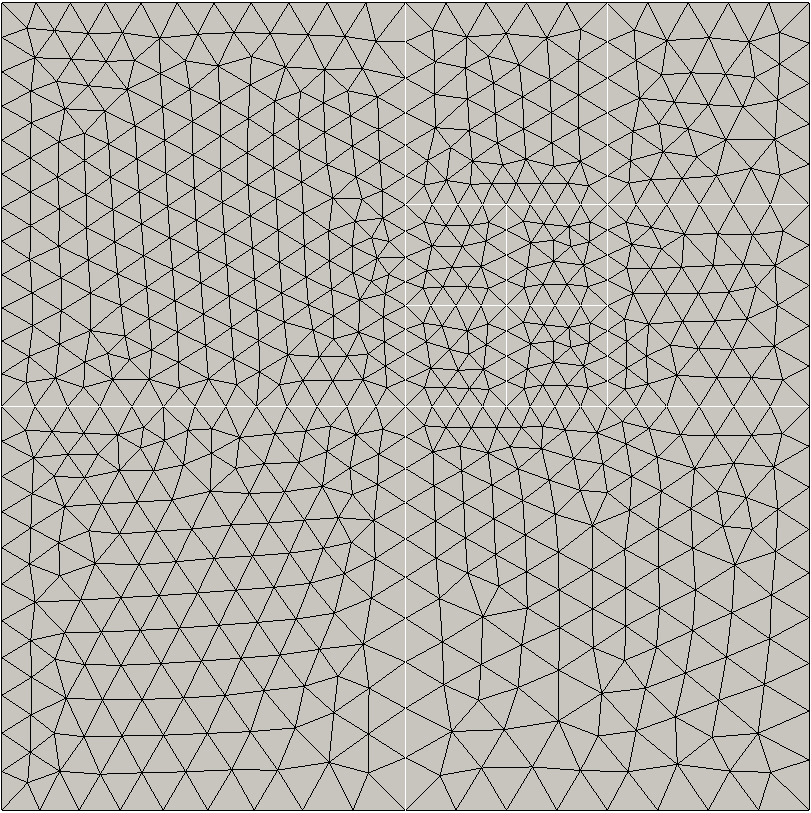}
}
\subfloat[EDFM]{
  \includegraphics[width=0.315\textwidth]{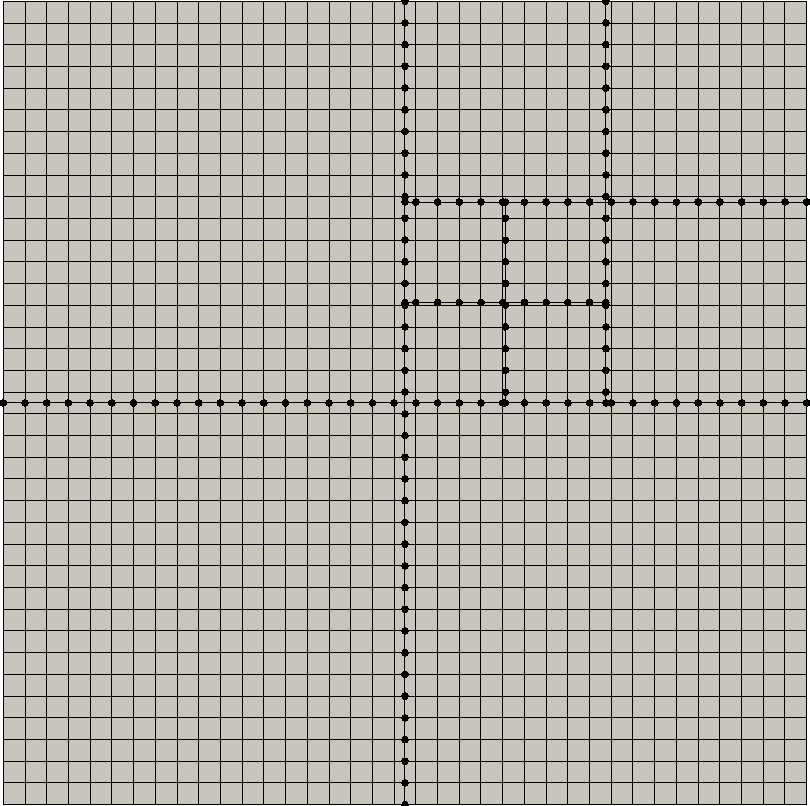}
}
\\
\subfloat[mortar-DFM]{
  \includegraphics[width=0.315\textwidth]{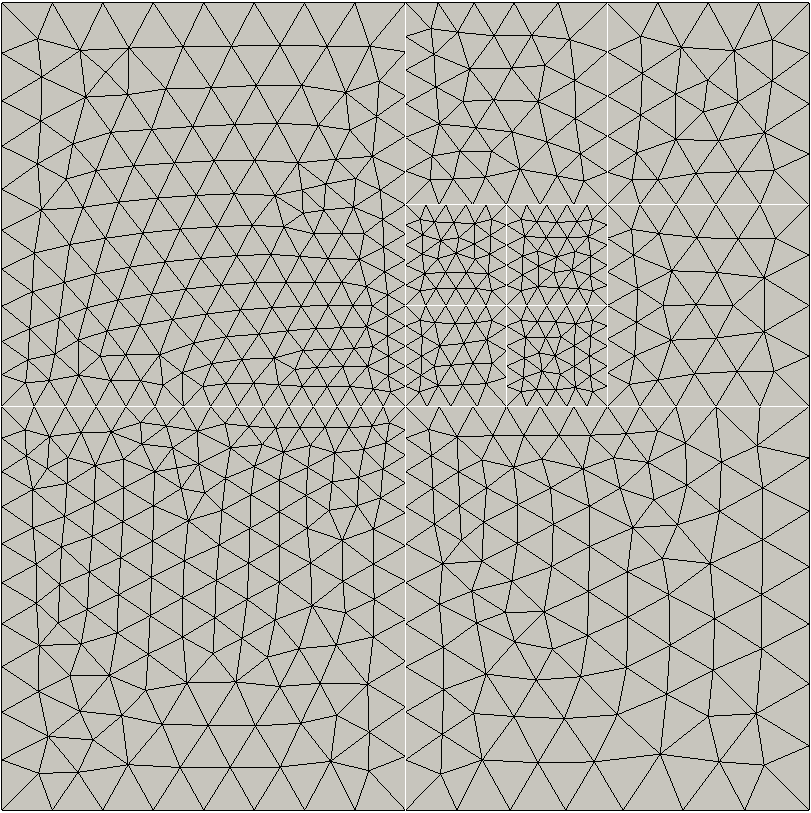}
}
\subfloat[P-XFEM]{
  \includegraphics[width=0.315\textwidth]{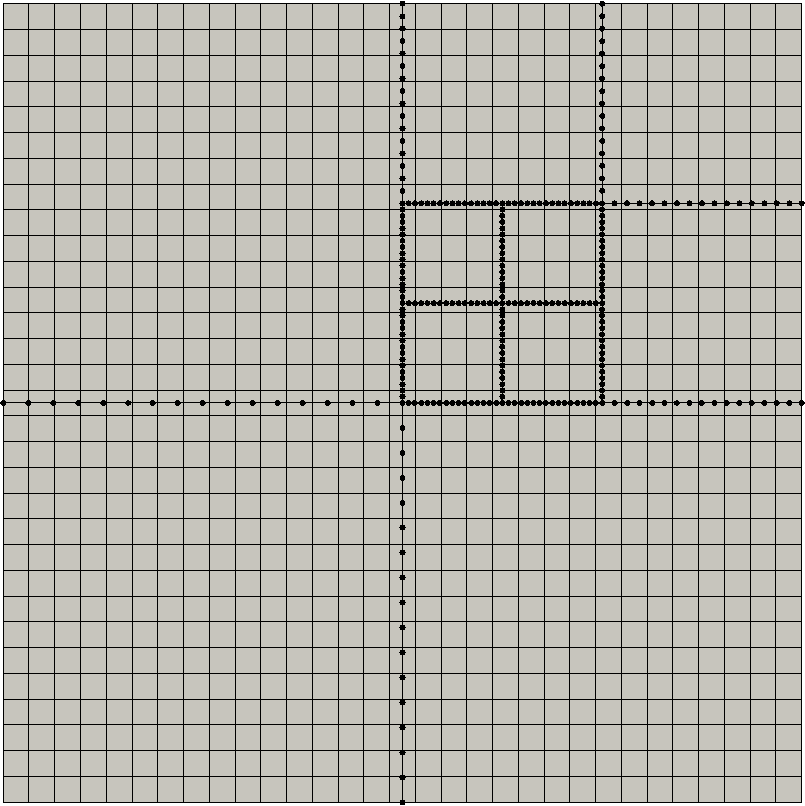}
}
\subfloat[D-XFEM]{
  \includegraphics[width=0.315\textwidth]{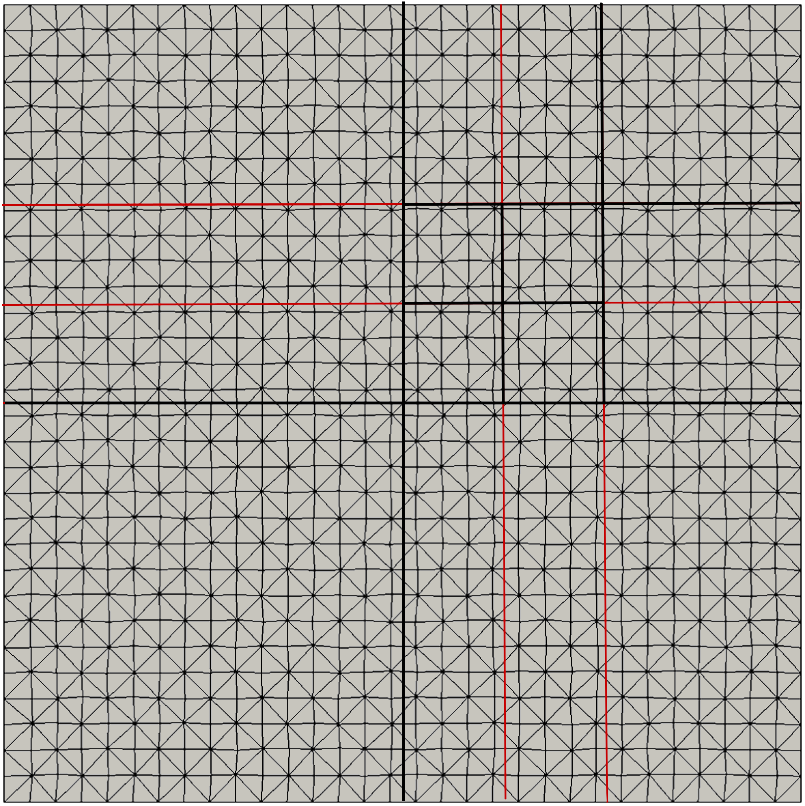}
}
\caption{Benchmark \theBenchmarkCounter: the grids used by the different methods. In the DXFEM grid the red lines indicate the virtual extension of the fractures up to the boundary.}
\label{fig:geiger_grids}
\end{figure}

The first distinction between the different schemes are given in Table \ref{tbl:geiger_grids}, where the number of degrees of freedom, matrix elements and fracture elements for all the participating methods are listed. The corresponding grids are visualized in Figure \ref{fig:geiger_grids}.
\renewcommand{\arraystretch}{1.1}
\begin{table}[hbt]
\centering
\begin{tabular}{|l|c|c|c|}\hline
\textbf{method} & \textbf{d.o.f.}
& \textbf{matrix elements} & \textbf{fracture elements} \\\hline
Box-DFM    & 1422 & 2691 triangles & 130 \\\hline
CC-DFM     & 1481 & 1386 triangles & 95\\\hline
EDFM       & 1501 & 1369 quads & 132 \\\hline
mortar-DFM & 3366 & 1280 triangles & 75\\\hline
P-XFEM     & 1632 & 961 quads & 318 \\\hline
D-XFEM     & 4474 & 1250 triangles & 126 \\\hline
MFD        & 2352280 & 1136456 quads & 38600 \\\hline
\end{tabular}
\caption{Grids for Benchmark \theBenchmarkCounter.}
\label{tbl:geiger_grids}
\end{table}
\renewcommand{\arraystretch}{1.0}

\subsubsection{Conductive Fracture Network}
\label{sec:geiger_example1}

\begin{figure}[hbt]
\centering
\includegraphics[width=0.5\linewidth]{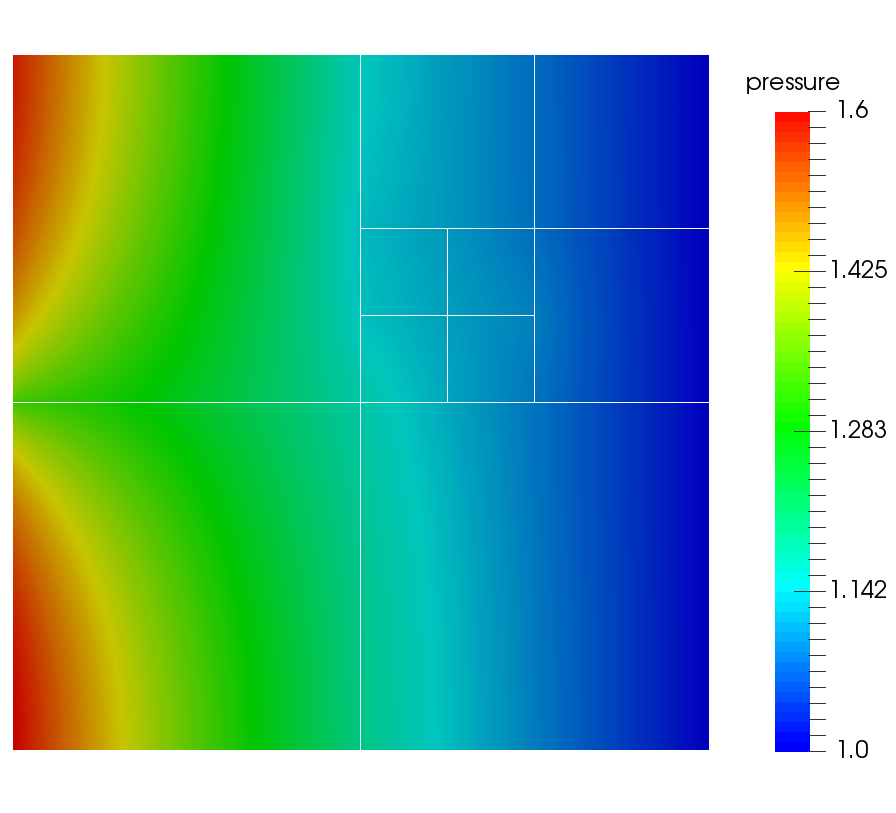}
\caption{Benchmark \theBenchmarkCounter\xspace with conductive fractures: pressure reference solution.}
\label{fig:geiger_conductive_reference}
\end{figure}

First, we consider a highly conductive network by setting $\kfn=\kft=10^4$.
The pressure distribution of the corresponding reference solution is shown in Figure \ref{fig:geiger_conductive_reference}.

\begin{figure}[hbt]
\centering
\subfloat[Horizontal line at $y=0.7$.]{
\includegraphics[width=0.49\textwidth]{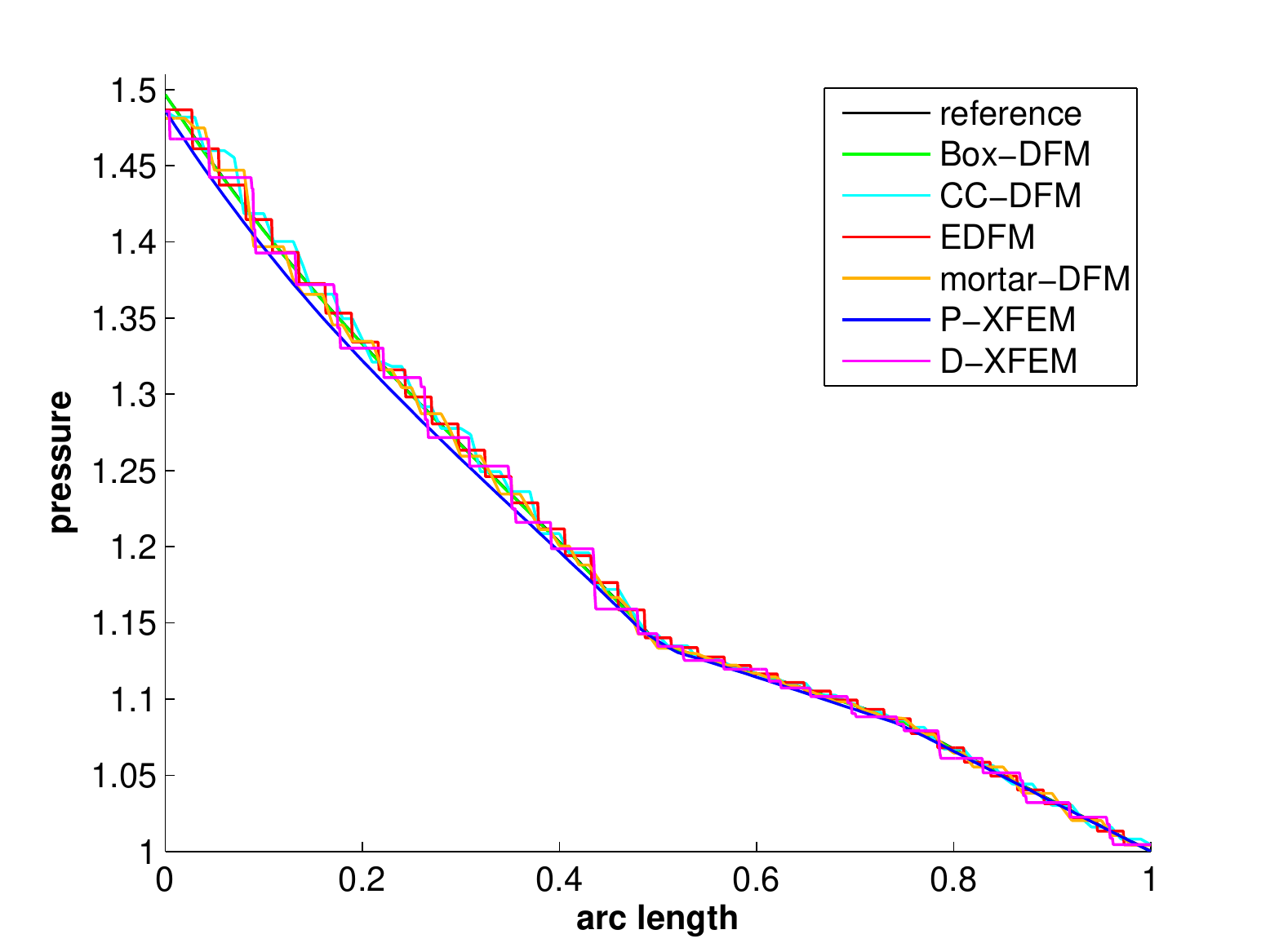}
}
\subfloat[Longest vertical fracture at $x=0.5$.]{
\includegraphics[width=0.49\textwidth]{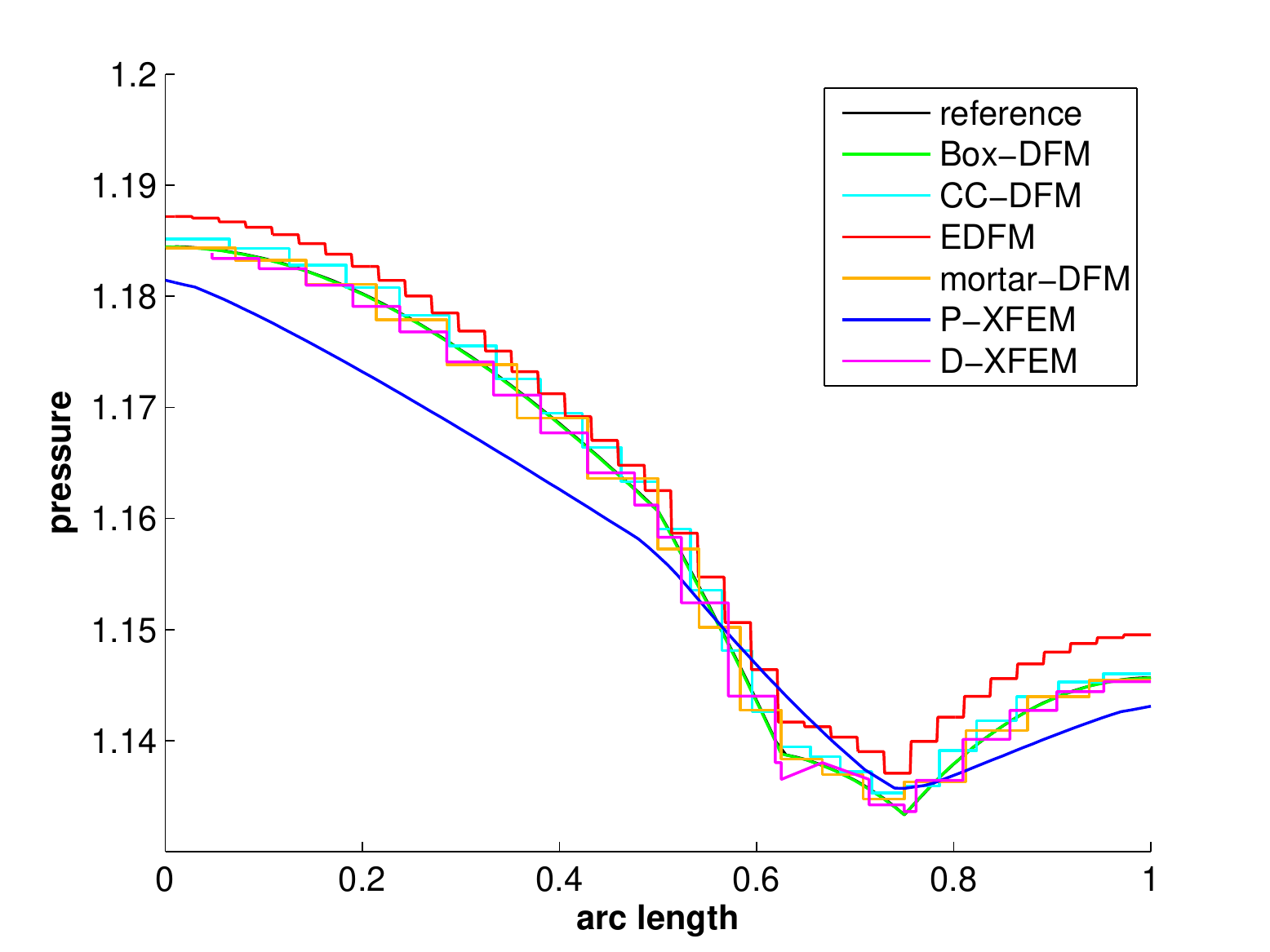}
}
\caption{Benchmark \theBenchmarkCounter\xspace with conductive fractures: comparison of values along two lines.}
\label{fig:permeable_lineplots}
\end{figure}

The pressure distributions given by the different methods are first compared along two lines, one horizontal at $y=0.7$ and one vertical at $x=0.5$. As shown in Figure \ref{fig:permeable_lineplots}, all results are relatively close to the reference solution. Qualitatively, we observe that P-XFEM produces a more diffuse pressure profile in the vertical fracture.

\renewcommand{\arraystretch}{1.1}
\begin{table}[hbt]
\centering
\begin{tabular}{|l|c|c|c|c|}\hline
\textbf{method} & $err_\text{m}$& $err_\text{f}$
& \textbf{nnz/size$^2$} & $\|\cdot\|_2$\textbf{-cond} \\\hline
Box-DFM    & 6.7e-3 & 1.1e-3 & 4.7e-3 & 7.9e3 \\\hline
CC-DFM     & 1.1e-2 & 5.0e-3 & 2.7e-3 & 5.6e4 \\\hline
EDFM       & 6.5e-3 & 4.0e-3 & 3.3e-3 & 5.6e4 \\\hline
mortar-DFM & 1.0e-2 & 7.4e-3 & 1.8e-3 & 2.4e6 \\\hline
P-XFEM     & 1.7e-2 & 6.0e-3 & 7.8e-3 & 6.8e9 \\\hline
D-XFEM     & 9.6e-3 & 8.9e-3 &  1.3e-3 & 1.2e6\\\hline
\end{tabular}
\caption{Discretization errors and matrix characteristics for Benchmark \theBenchmarkCounter\xspace with conductive fractures.}
\label{tbl:permeable_error}
\end{table}
\renewcommand{\arraystretch}{1.0}

The results for the conducting fractures are similar to those presented in the
first benchmark. In particular, the performance of the methods is comparable as
shown by both the matrix and the fracture errors. In fact, since the degree of
sparsity does not differ significantly either, the only notable differences
between the methods are the number of degrees of freedom and the condition
numbers, as shown in Table \ref{tbl:permeable_error}. In that context, the
mortar-DFM and D-XFEM are the clear outliers, containing a large number of
degrees of freedom due to the incorporated flux variable and resulting in high
condition numbers. The P-XFEM scheme exhibits the highest condition number, yet
we emphasize that it has significantly fewer degrees of freedom.

\subsubsection{Blocking Fracture Network}
\label{sec:geiger_example2}

\begin{figure}[hbt]
\centering
\includegraphics[width=0.5\linewidth]{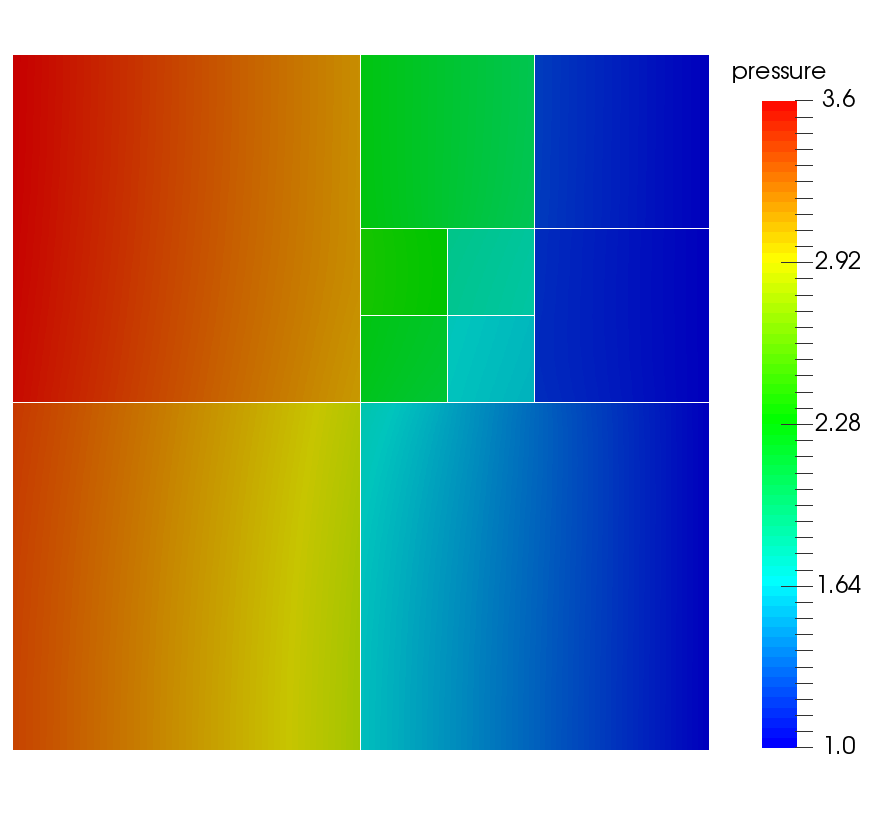}
\caption{Benchmark \theBenchmarkCounter\xspace with blocking fractures: pressure reference solution.}
\label{fig:geiger_blocking_reference}
\end{figure}

We now assume a blocking fracture network by setting $\kfn=\kft=10^{-4}$.
The pressure distribution of the corresponding reference solution is shown in Figure \ref{fig:geiger_blocking_reference}. The results clearly show the pressure discontinuities reminiscent of the low fracture permeability.

\begin{figure}[hbt]
\centering
\includegraphics[width=0.9\textwidth]{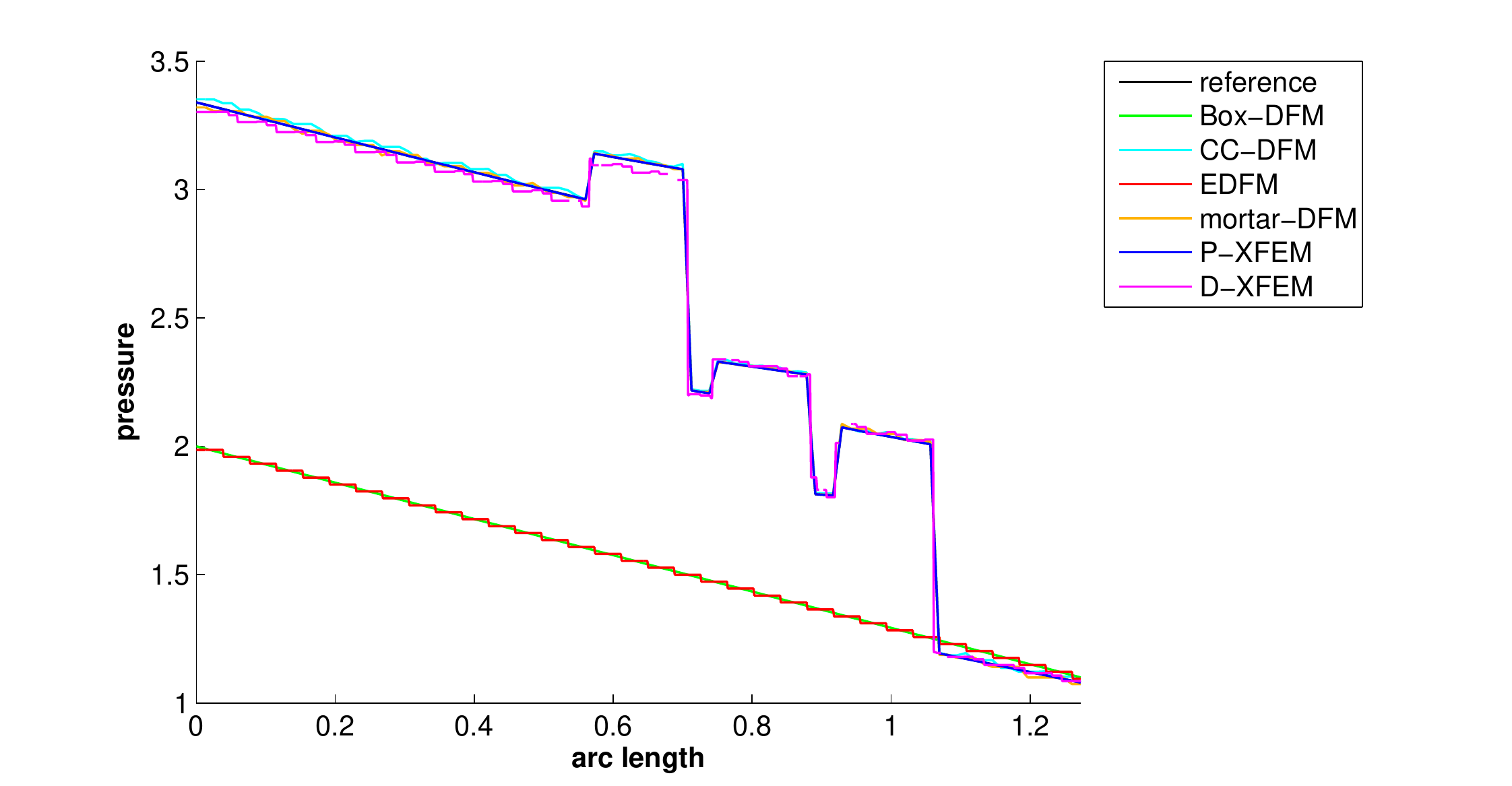}
\caption{Benchmark \theBenchmarkCounter\xspace with blocking fractures: values along the line $(0.0,0.1)-(0.9,1.0)$.}
\label{fig:geiger_imp}
\end{figure}

Figure \ref{fig:geiger_imp} compares the results of the different methods along a
diagonal line crossing the whole domain from $(0.0,0.1)$ to $(0.9,1.0)$. The discretization errors, sparsity densities, and condition numbers for the different
methods are given in Table \ref{tbl:blocking_error}.

\renewcommand{\arraystretch}{1.1}
\begin{table}[hbt]
\centering
\begin{tabular}{|l|c|c|c|c|}\hline
\textbf{method} & $err_\text{m}$ & $err_\text{f}$
& \textbf{nnz/size$^2$} & $\|\cdot\|_2$\textbf{-cond} \\\hline
Box-DFM    & 4.1e-1 & 3.8e-1 & 4.7e-3 & 3.5e3 \\\hline
CC-DFM     & 5.7e-3 & 4.6e-3 & 2.7e-3 & 2.6e4 \\\hline
EDFM       & 2.9e-1 & 3.2e-1 & 3.3e-3 & 9.2e3 \\\hline
mortar-DFM & 4.5e-3 & 4.9e-3 & 1.6e-3 & 9.0e2 \\\hline
P-XFEM     & 2.9e-3 & 2.2e-2 & 8.1e-3 & 2.0e4 \\\hline
D-XFEM     & 1.0e-2 & 1.9e-2  & 1.3e-3 & 2.2e6 \\\hline
\end{tabular}
\caption{Discretization errors and matrix characteristics for Benchmark \theBenchmarkCounter\xspace with blocking fractures.}
\label{tbl:blocking_error}
\end{table}
\renewcommand{\arraystretch}{1.0}

In the case of blocking fractures, the distinction between the different methods is more apparent. As mentioned above, the Box-DFM and EDFM schemes are unable to capture the resulting pressure discontinuities. As a result, these methods show large errors in both the matrix and the fracture domains. The remaining methods, which are capable of handling discontinuities, do not differ significantly among each other in terms of fracture and matrix errors. We do note that the condition numbers have improved significantly for the mortar-DFM and P-XFEM schemes. Conversely, for CC-DFM and D-XFEM, condition numbers for the blocking fractures case are similar to those obtained for the permeable fractures case.

\clearpage
\stepcounter{BenchmarkCounter}
\subsection{Benchmark \theBenchmarkCounter: Complex Fracture Network}
\label{sec:anna}

\begin{figure}[hbt]
    \centering
    \subfloat[]{
        \centering
        \resizebox{0.4\textwidth}{!}{\fontsize{20pt}{7.2}\selectfont%
        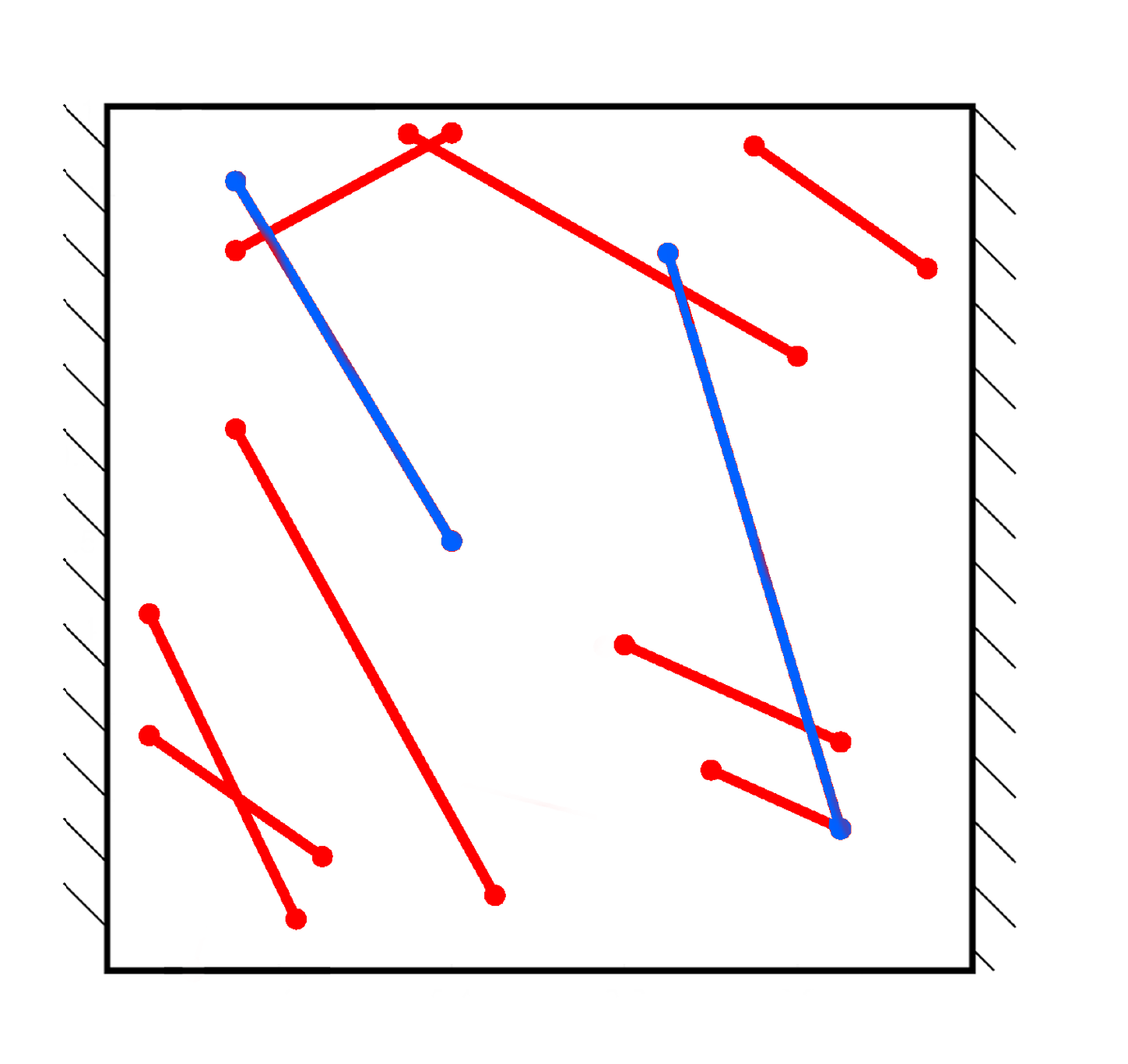}}%
    \hspace*{0.1\textwidth}%
    \subfloat[]{
        \centering
        \resizebox{0.4\textwidth}{!}{\fontsize{20pt}{7.2}\selectfont%
        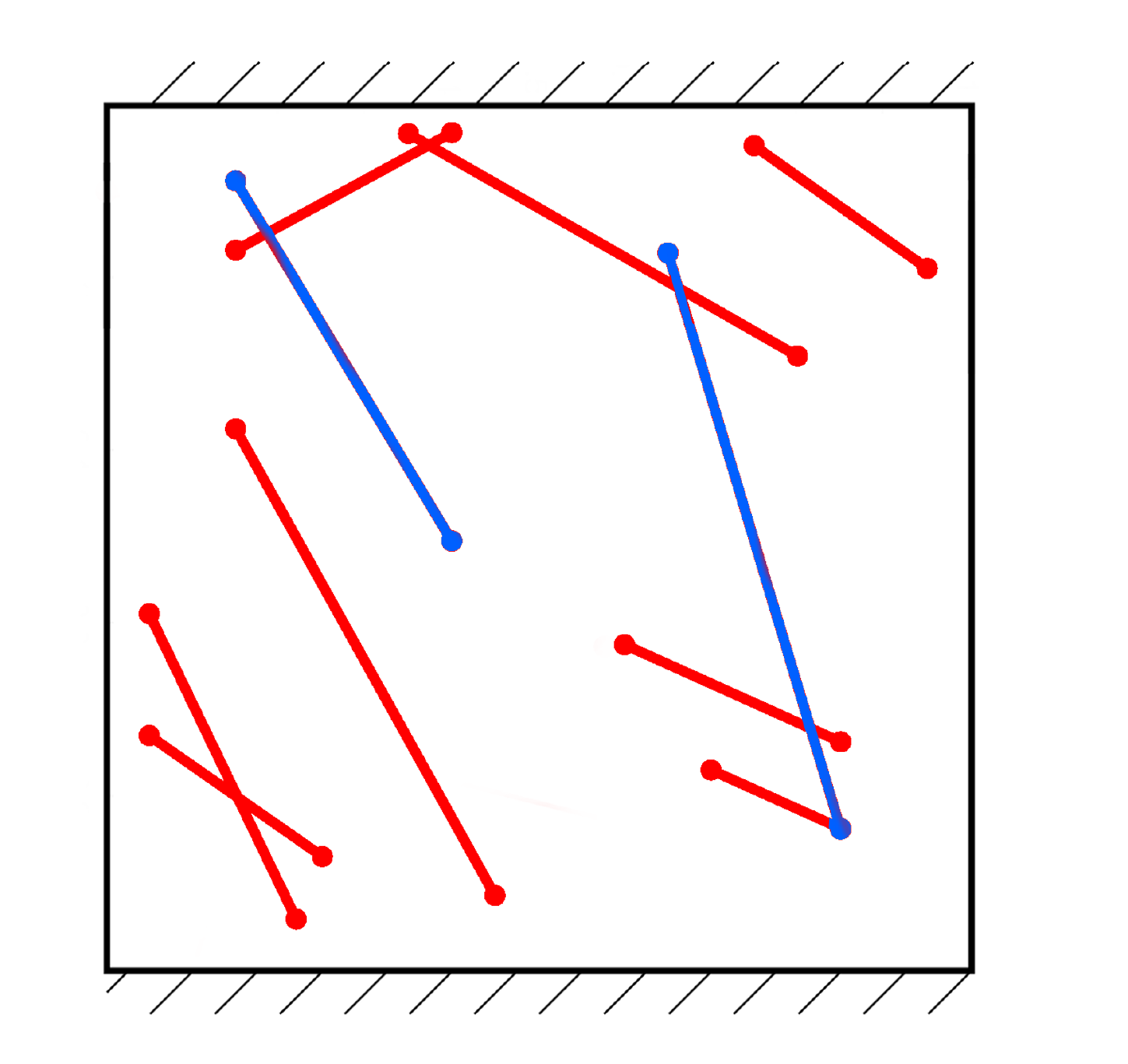}
    }%
\caption{Benchmark \theBenchmarkCounter: Domain and boundary conditions for
cases (a) and (b). The red fractures are conductive, the blue ones are blocking.}%
\label{fig:complex_domain_anna}
\end{figure}

\begin{figure}[hbt]
\subfloat[]{
  \includegraphics[width=0.45\textwidth]{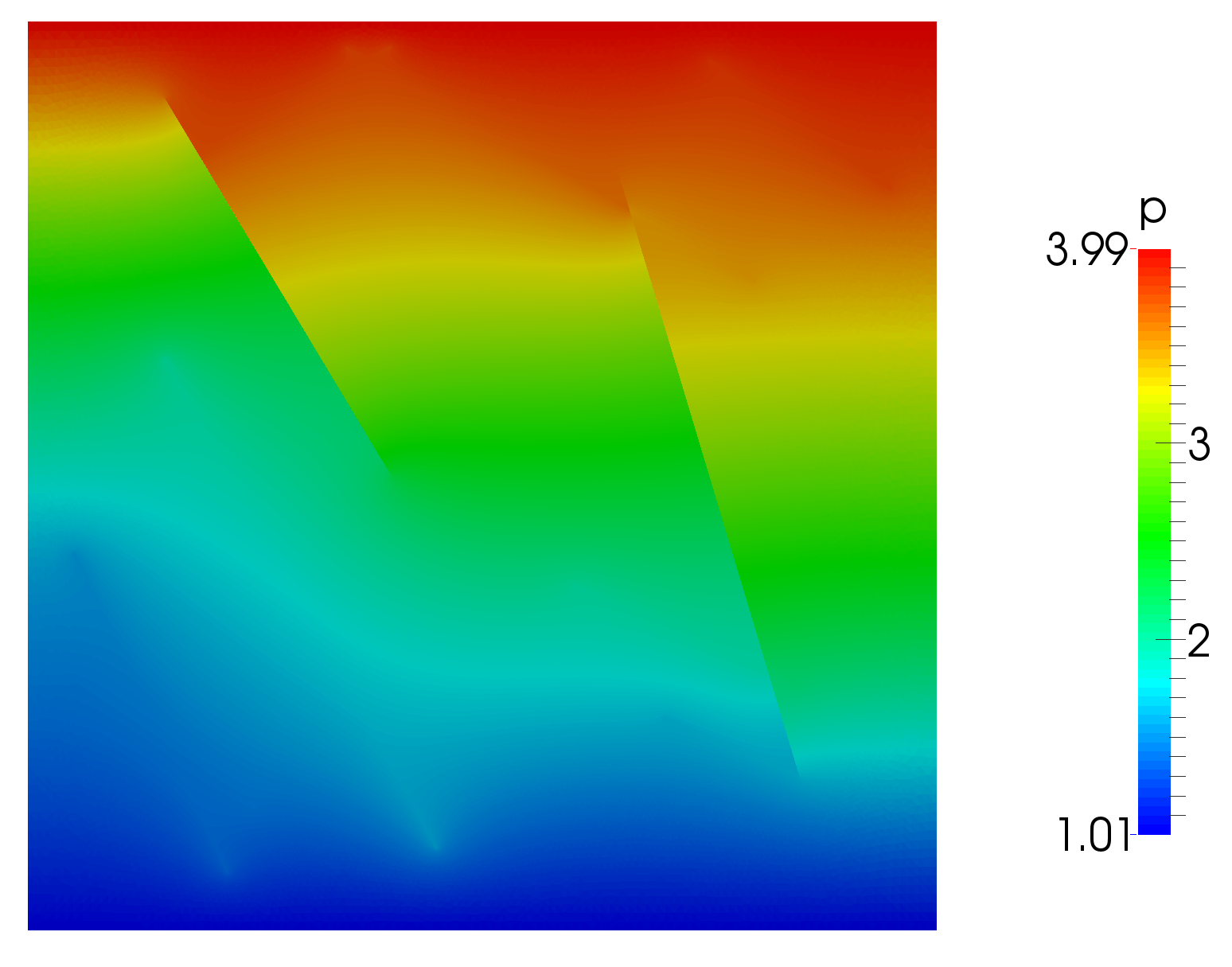}
}
\subfloat[]{
  \includegraphics[width=0.45\textwidth]{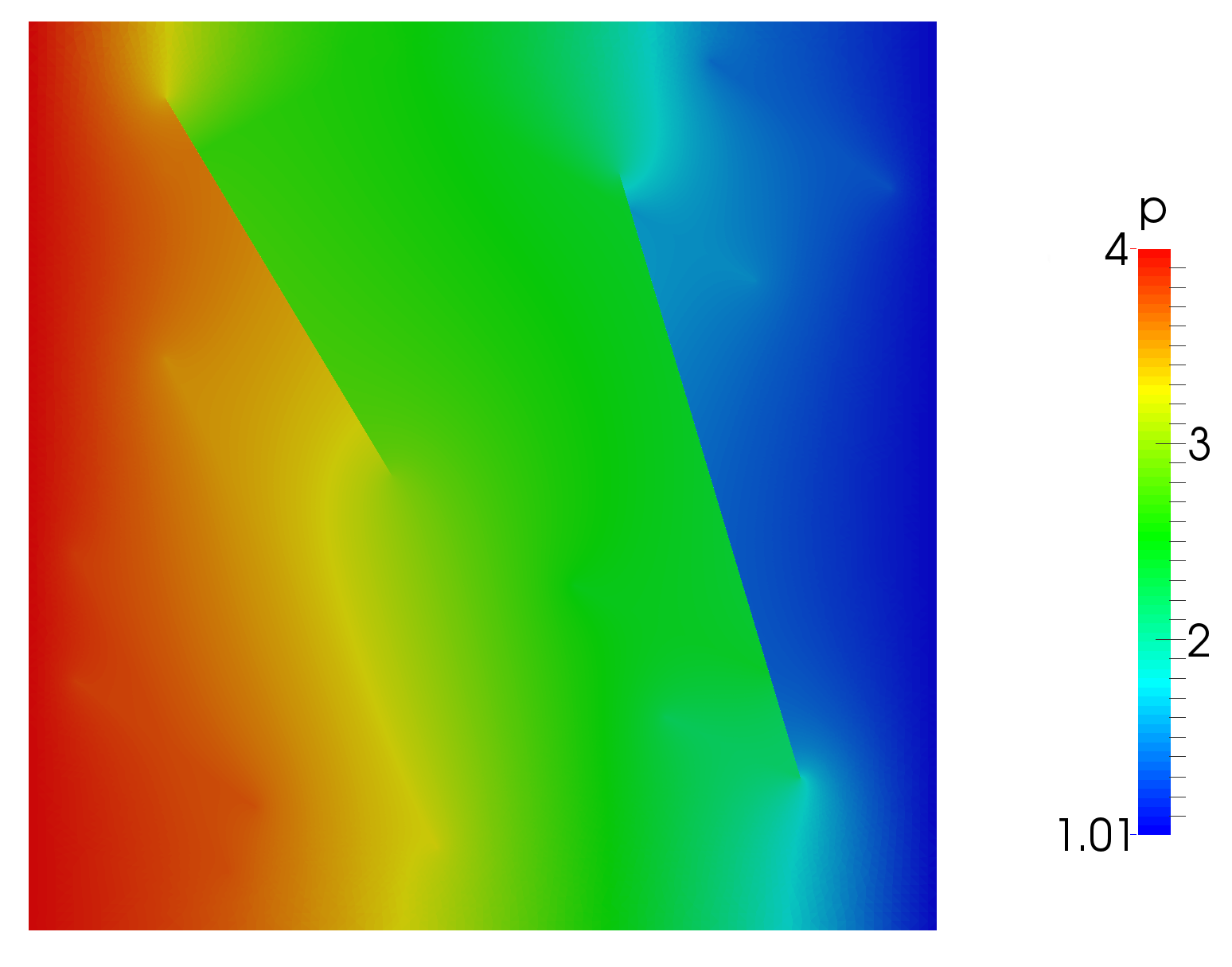}
}
\caption{Benchmark \theBenchmarkCounter: reference solution for cases a) and b)}
\label{fig:refsol_anna}
\end{figure}

This test case considers a small but complex fracture network that includes permeable and blocking fractures. The domain and boundary conditions are shown in Figure \ref{fig:complex_domain_anna}.
The exact coordinates for the fracture positions are provided in \ref{sec:anna_coord}. The fracture network contains ten straight immersed fractures, grouped in disconnected networks. The aperture is $\aperture=10^{-4}$ for all fractures, and permeability is $\kfn=\kft=10^4$ for all fractures except for fractures $4$ and $5$ which are blocking fractures with $\kfn=\kft=10^{-4}$. Note that we are considering two subcases a) and b) with a pressure gradient which is predominantly vertical and horizontal respectively, to better highlight the impact of the blocking fractures.
The corresponding reference solutions are depicted in Figure \ref{fig:refsol_anna}.

\begin{figure}[hbt]
\centering
\subfloat[Box-DFM]{
  \includegraphics[width=0.315\textwidth]{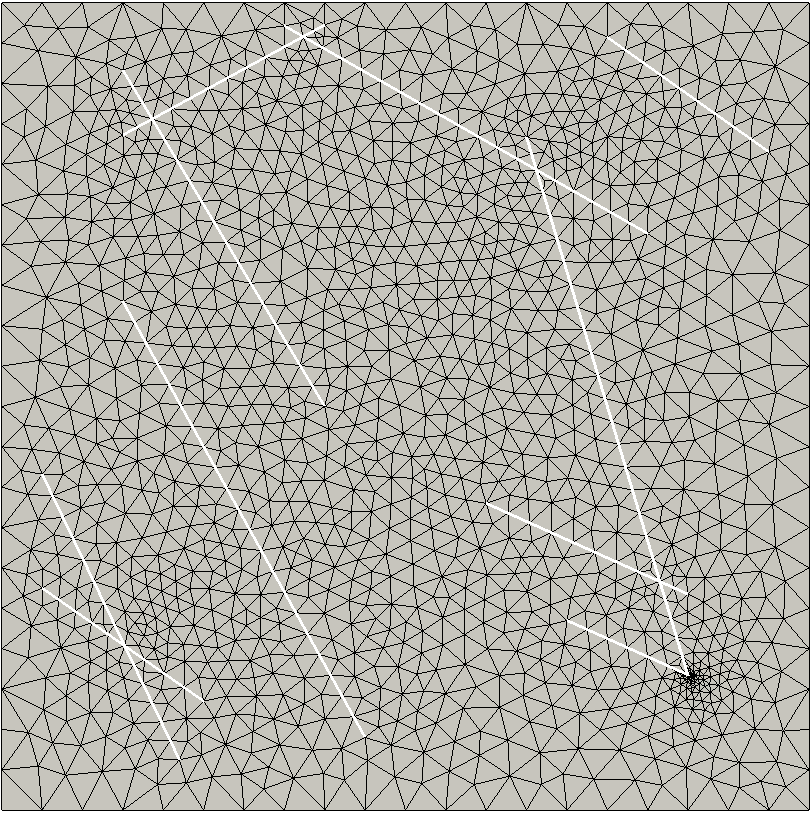}
}
\subfloat[CC-DFM]{
  \includegraphics[width=0.315\textwidth]{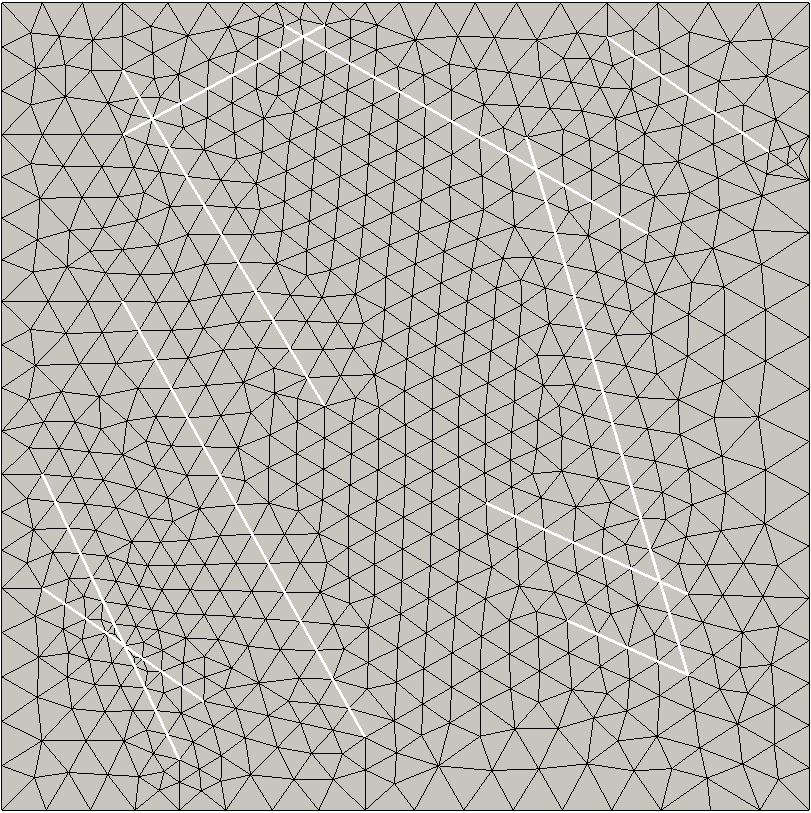}
}
\subfloat[EDFM]{
  \includegraphics[width=0.315\textwidth]{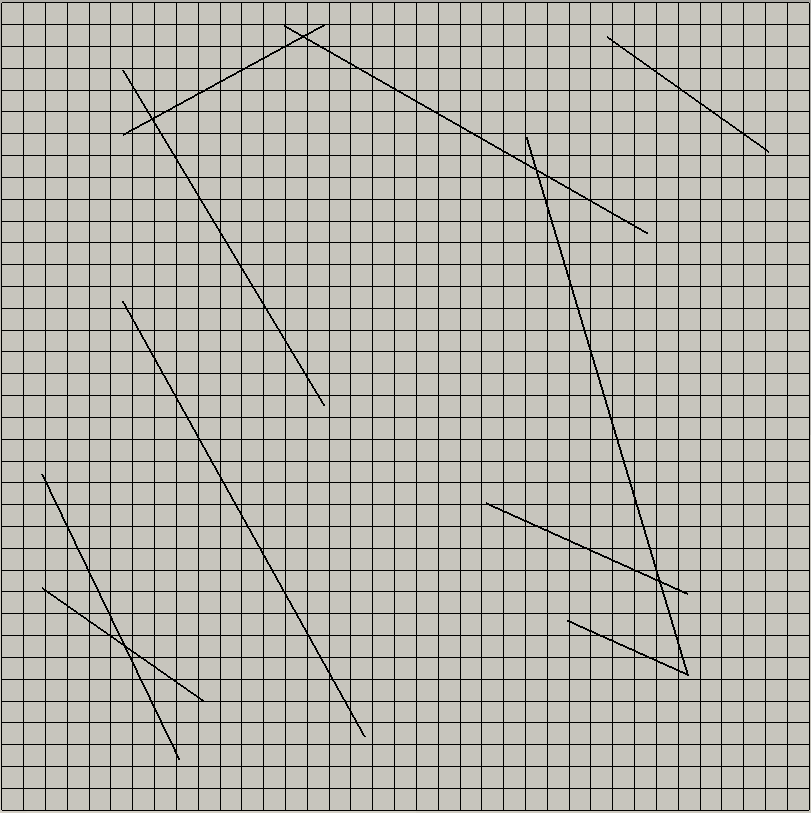}
}
\\
\subfloat[mortar-DFM]{
  \includegraphics[width=0.315\textwidth]{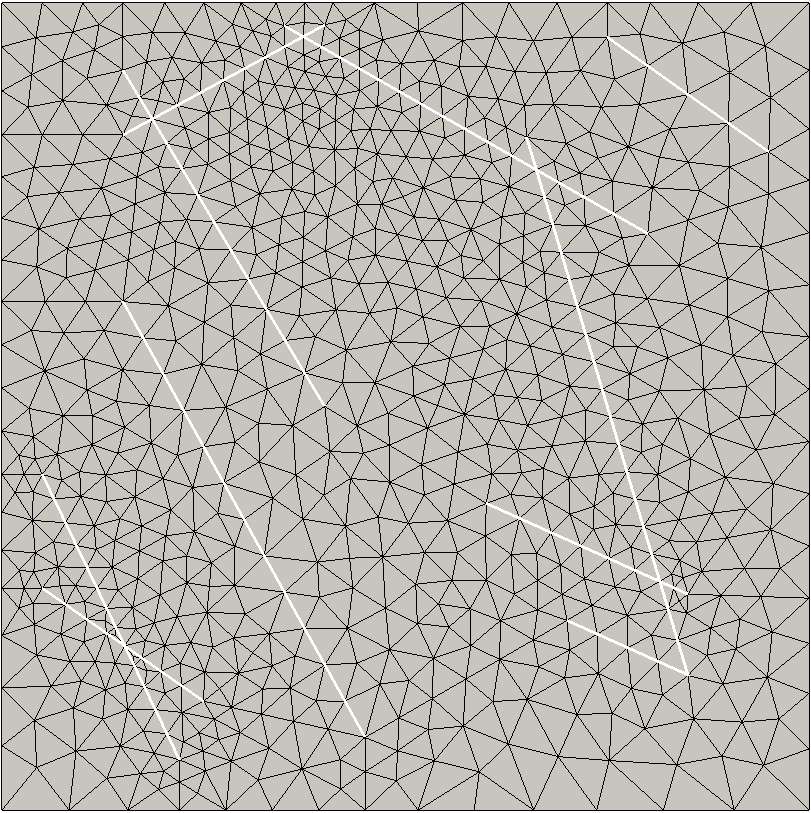}
}
\subfloat[P-XFEM: n/a]{
  \includegraphics[width=0.315\textwidth]{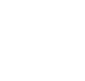}
}
\subfloat[D-XFEM]{
  \includegraphics[width=0.315\textwidth]{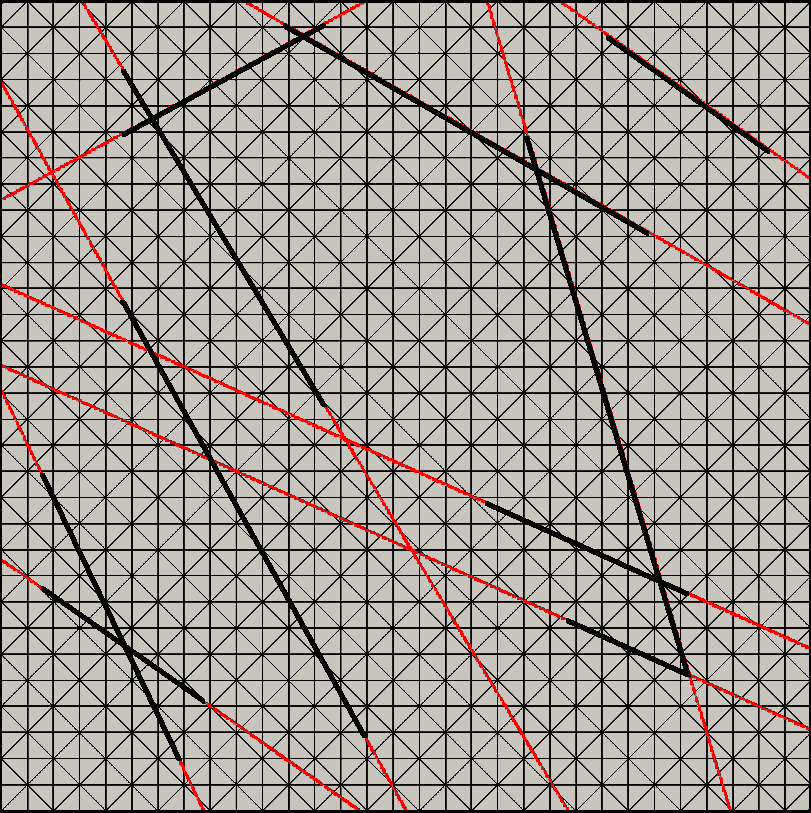}
}
\caption{Benchmark \theBenchmarkCounter: the grids used by the different methods. In the DXFEM grid the red lines indicate the virtual extension of the fractures up to the boundary.}
\label{fig:anna_grids}
\end{figure}

Table \ref{tbl:anna_grids} lists the number of degrees of freedom, matrix elements and fracture elements for all the participating methods.
The corresponding grids are visualized in Figure \ref{fig:anna_grids}.
\renewcommand{\arraystretch}{1.1}
\begin{table}[hbt]
\centering
\begin{tabular}{|l|c|c|c|}\hline
\textbf{method} & \textbf{d.o.f.}
& \textbf{matrix elements} & \textbf{fracture elements} \\\hline
Box-DFM    & 1460 & 2838 triangles & 155 \\\hline
CC-DFM     & 1510 & 1407 triangles & 103\\\hline
EDFM       & 1572 & 1369 quads & 203 \\\hline
mortar-DFM & 3953 & 1452 triangles & 105\\\hline
D-XFEM     & 7180 & 1922 triangles & 199 \\\hline
MFD        & 1800770 & 1192504 mixed & 7876 \\\hline
\end{tabular}
\caption{Grids for Benchmark \theBenchmarkCounter.}
\label{tbl:anna_grids}
\end{table}
\renewcommand{\arraystretch}{1.0}

The P-XFEM method could not participate in this benchmark example.
Its current implementation requires that each matrix element face is cut by
at most one fracture branch. While it would be possible to construct a matrix grid that
satisfies this requirement, this would contradict the promised advantage of admitting
independent fracture and matrix grids.

\subsubsection{Flow from Top to Bottom}

Table \ref{tbl:anna_error_toptobottom} lists the discretization errors for the first variant, namely, the flow from top to bottom.
\renewcommand{\arraystretch}{1.1}
\begin{table}[hbt]
\centering
\begin{tabular}{|l|c|c|c|c|}\hline
\textbf{method} & $err_\text{m}$ & $err_\text{f}$
& \textbf{nnz/size$^2$} & $\|\cdot\|_2$\textbf{-cond} \\\hline
Box-DFM    & 4.4e-2 & 3.8e-2 & 4.6e-3 & 4.5e3 \\\hline
CC-DFM     & 2.6e-2 & 3.3e-2 & 2.7e-3 & 3.8e4\\\hline
EDFM       & 3.8e-2 & 4.5e-2 & 3.1e-3 & 1.2e6\\\hline
mortar-DFM & 1.0e-2 & 1.7e-2 & 1.4e-3 & 1.1e6\\\hline
D-XFEM     & 1.9e-2 & 2.9e-2 & 8.2e-4 & 8.1e3\\\hline
\end{tabular}
\caption{Discretization errors and matrix characteristics for the first variant of Benchmark \theBenchmarkCounter.}
\label{tbl:anna_error_toptobottom}
\end{table}
\renewcommand{\arraystretch}{1.0}

Even though this is still a synthetic case, we can see that the geometry of the
network starts to be an issue: relatively small intersection angles are present,
for instance, between fractures 1 and 2. Another difficulty consists in the
coexistence of permeable and blocking fractures which intersect each other: on
one hand, some of the methods are not well suited to describe a blocking
behavior, on the other hand the coupling conditions at the intersection become
less trivial in these cases. All the participating methods that account
explicitly for the effect of permeability at the fracture intersections have
adopted the harmonic average in the case of a permeable and a blocking fracture
crossing each other. The errors reported in Table
\ref{tbl:anna_error_toptobottom} show that the methods requiring the continuity
of pressure (EDFM and the Box-DFM) exhibit slightly higher errors in the matrix.
However, the difference is not particularly sharp, since in this sub-case the
average pressure gradient is almost parallel to the blocking fractures.

\subsubsection{Flow from Left to Right}

The discretization errors for the second variant, namely, the flow from left to right, are summarized in Table \ref{tbl:anna_error_lefttoright}.
\renewcommand{\arraystretch}{1.1}
\begin{table}[hbt]
\centering
\begin{tabular}{|l|c|c|c|c|}\hline
\textbf{method} & $err_\text{m}$ & $err_\text{f}$
& \textbf{nnz/size$^2$} & $\|\cdot\|_2$\textbf{-cond} \\\hline
Box-DFM    & 7.5e-2 & 7.0e-2 & 4.6e-3 & 5.6e3 \\\hline
CC-DFM     & 5.2e-2 & 7.3e-2 & 2.7e-3 & 4.5e4 \\\hline
CC-DFM*     & 1.1e-2 & 2.7e-2 & 2.6e-3 & 8.1e5 \\\hline
EDFM       & 5.8e-2 & 8.9e-2 & 3.1e-3 & 1.2e6 \\\hline
mortar-DFM & 1.3e-2 & 2.7e-2 & 1.4e-3 & 7.3e8 \\\hline
D-XFEM     & 2.2e-2 & 3.6e-2 & 8.2e-4 & 8.1e3\\\hline
\end{tabular}
\caption{Discretization errors and matrix characteristics for the second variant of Benchmark \theBenchmarkCounter.}
\label{tbl:anna_error_lefttoright}
\end{table}
\renewcommand{\arraystretch}{1.0}

In this second case, since we impose pressure on the sides of the square domain, the solution is more challenging as we can observe from Figure \ref{fig:refsol_anna} and the gap between continuous and discontinuous methods increases. However, it should be noted that the errors remain of the same order of magnitude, indicating that all the methods capture the overall trend of the solution. The elimination of the fracture intersection cells in the CC-DFM is ill-suited for cases where fractures of different permeability cross. Therefore, we include a solution CC-DFM* in which we have not performed the removal for case b). The new results, reported in  Table \ref{tbl:anna_error_lefttoright},  show a far smaller error compared to the CC-DFM with elimination, but also demonstrate that the elimination significantly reduces the condition number.

\clearpage
\clearpage
\stepcounter{BenchmarkCounter}
\subsection{Benchmark \theBenchmarkCounter: a Realistic Case}
\label{sec:realistic}

\begin{figure}
    \centering
    \includegraphics[width=0.5\textwidth]{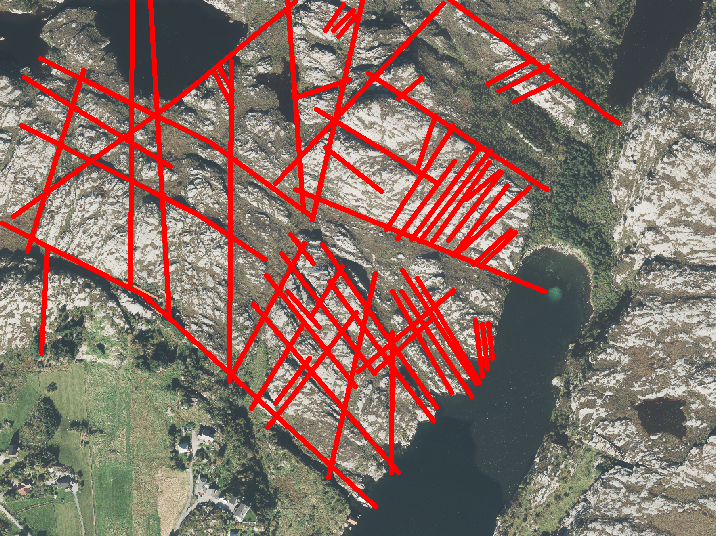}%
    \hfill%
    \includegraphics[width=0.45\textwidth]{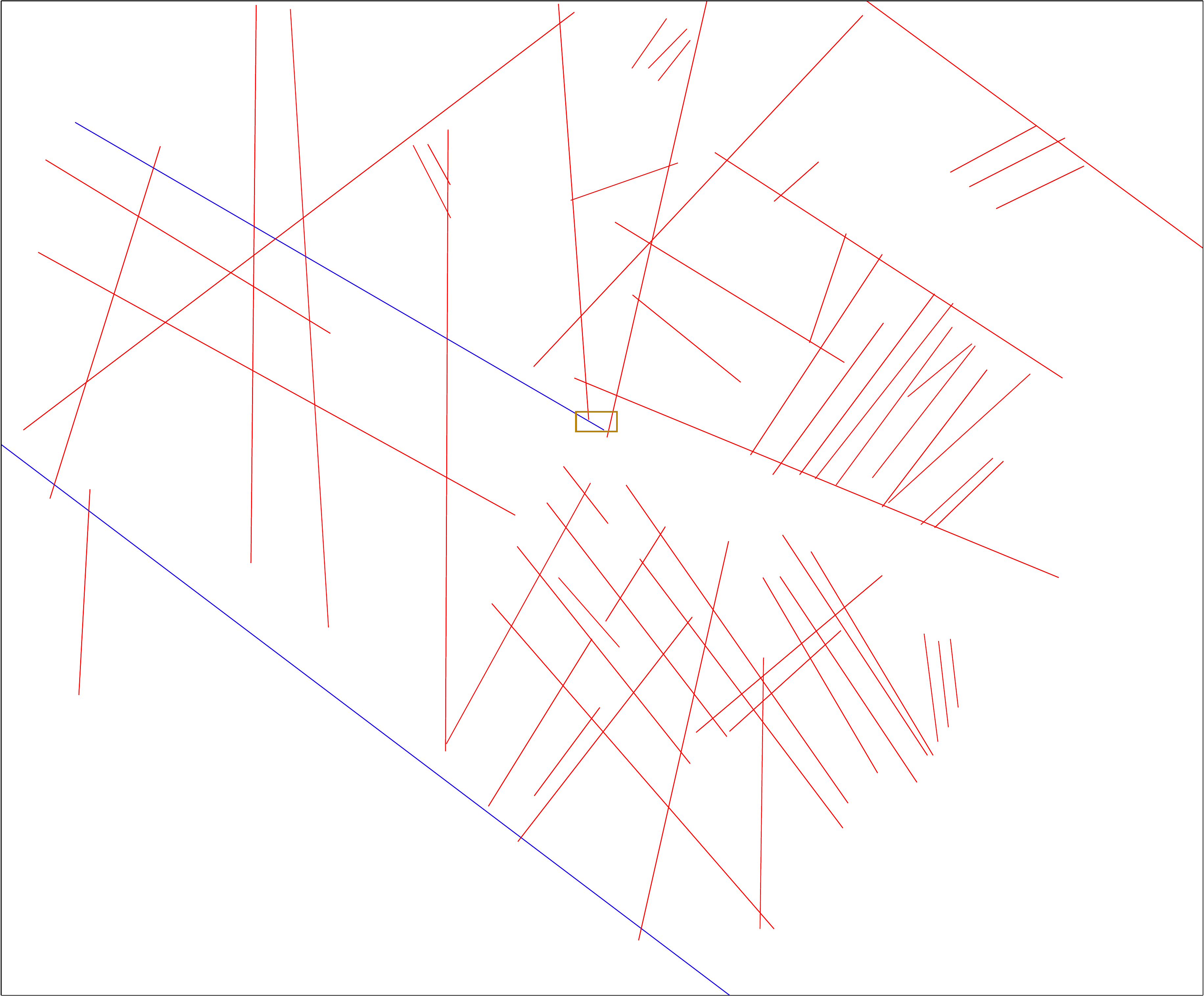}
    \caption{In the left the interpretation of the set of fractures superimposed
    to the map. In the right the geometry used in the simulations. The
    rectified fractures are depicted in blue.}%
    \label{fig:sotra}
\end{figure}

In this last test case we consider a real set of fractures from an
interpreted outcrop in the Sotra island, near Bergen in Norway. The set
is composed of 64 fractures grouped in 13 different connected networks, ranging
from isolated fractures up to tens of fractures each. In the interpretation process
two fractures were composed by more than one segment. However, since the
implementation of some methods relay on the fact that one fracture is
represented by a single geometrical object, we substitute them by a single segment.
It is worth to notice that we are
changing the connectivity of the system, nevertheless our goal is to make a
comparison of the previous schemes on a complex set of fractures. The
interpreted outcrop and the corresponding set of fractures are
represented in Figure \ref{fig:sotra}.
The size of the domain is $\unit[700]{m} \times \unit[600]{m}$ with uniform scalar permeability
equal to $\unit[10^{-14}]{m^2}$. For simplicity all the fractures have the same scalar
permeability equal in the tangential and normal direction to $\unit[10^{-8}]{m^2}$, and
aperture $\unit[10^{-2}]{m}$. We consider no-flow boundary condition on top and bottom,
pressure $\unit[1013250]{Pa}$ on the left, and pressure $\unit[0]{Pa}$ on the right of the
boundary of the domain. Due to the high geometrical complexity of the fracture
network not all involved numerical schemes/simulators could be used.
Nevertheless, it is worth to point out that for the others the main difficulty
in handling such geometry is an implementation issue rather than a limitation of
the scheme.

\begin{figure}
    \centering
    \subfloat[Box-DFM]{
    \includegraphics[width=0.5\textwidth]{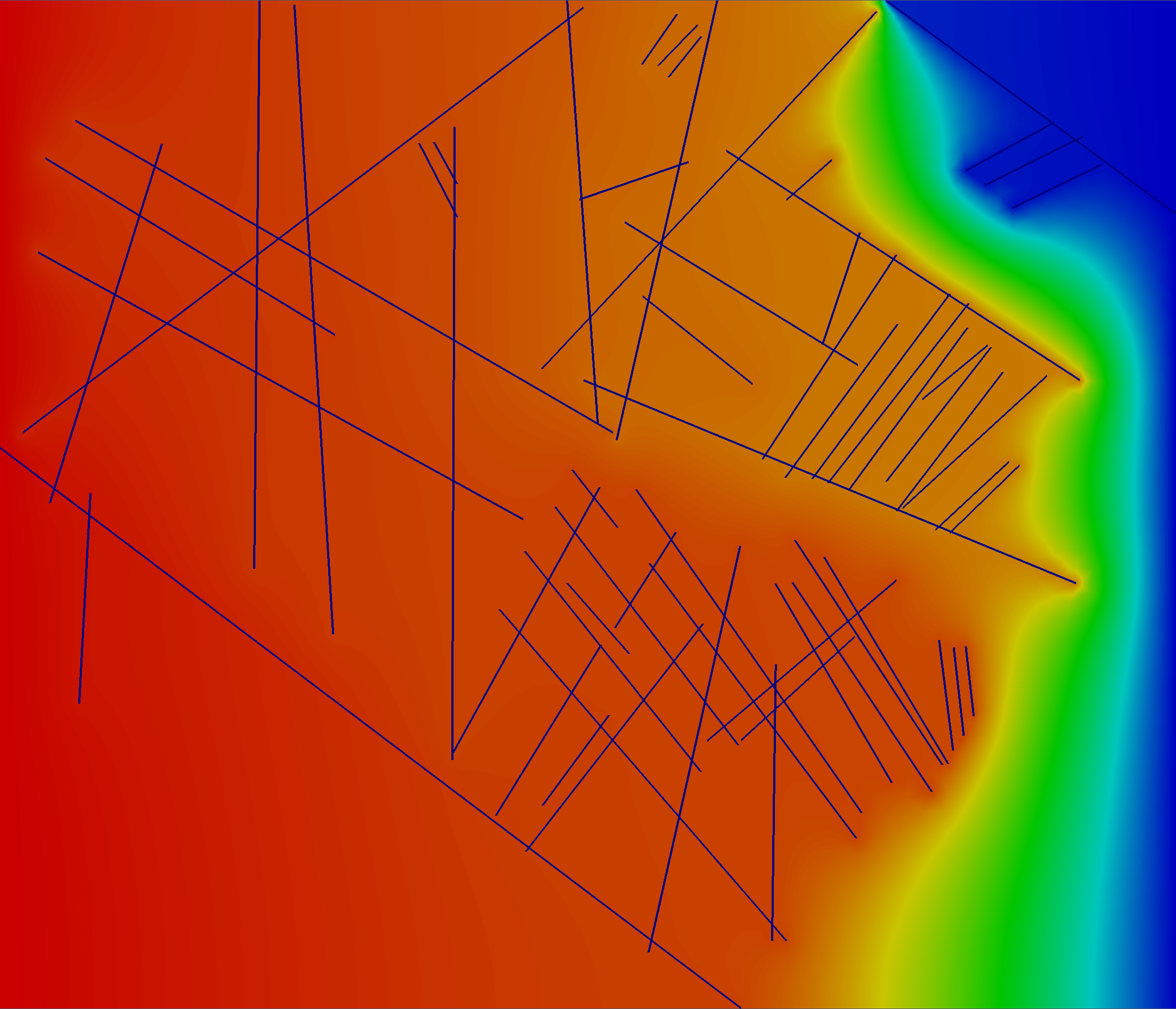}%
    }
    \subfloat[CC-DFM]{
    \includegraphics[width=0.5\textwidth]{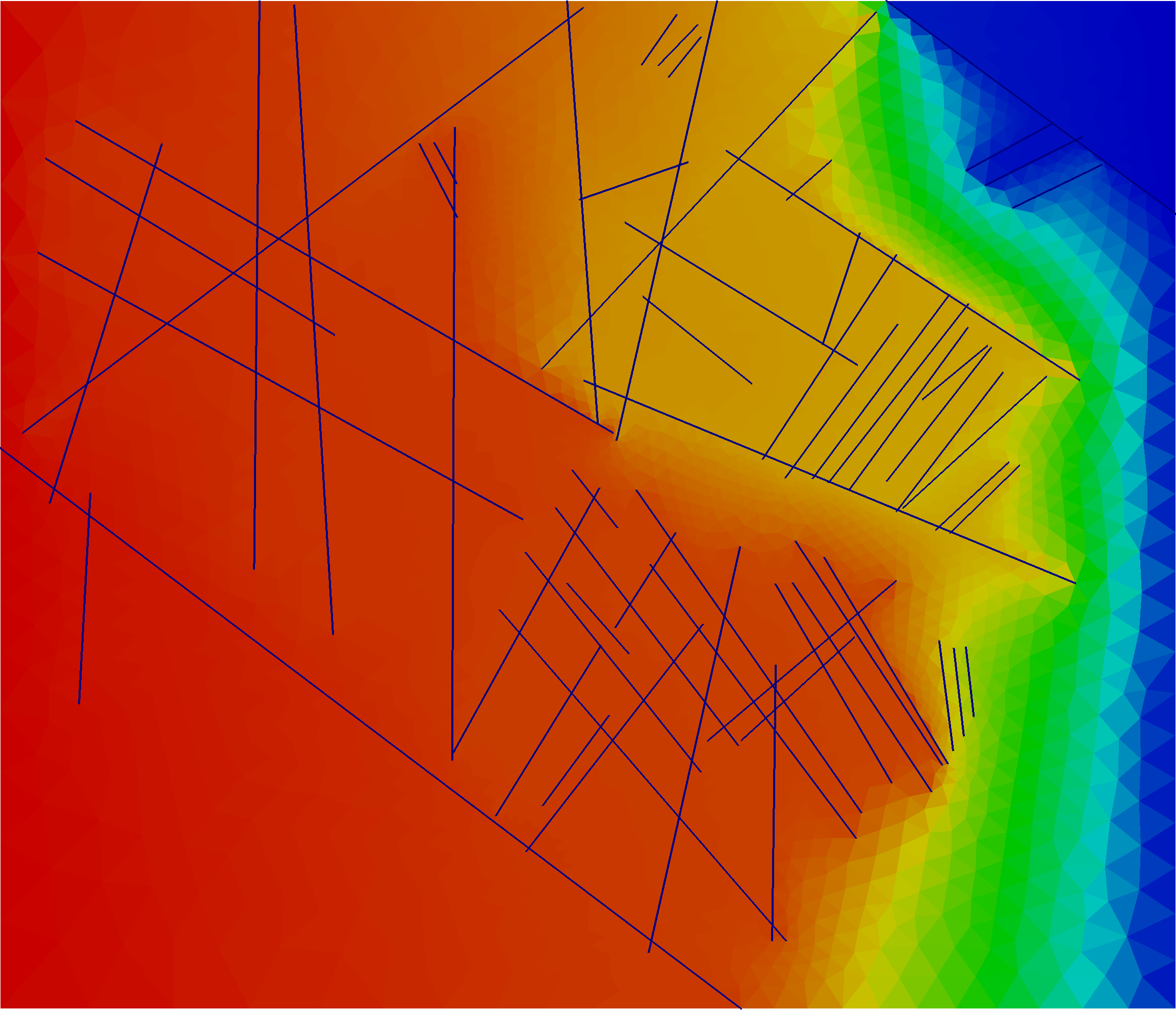}
    }\\
    \subfloat[EDFM]{
    \includegraphics[width=0.5\textwidth]{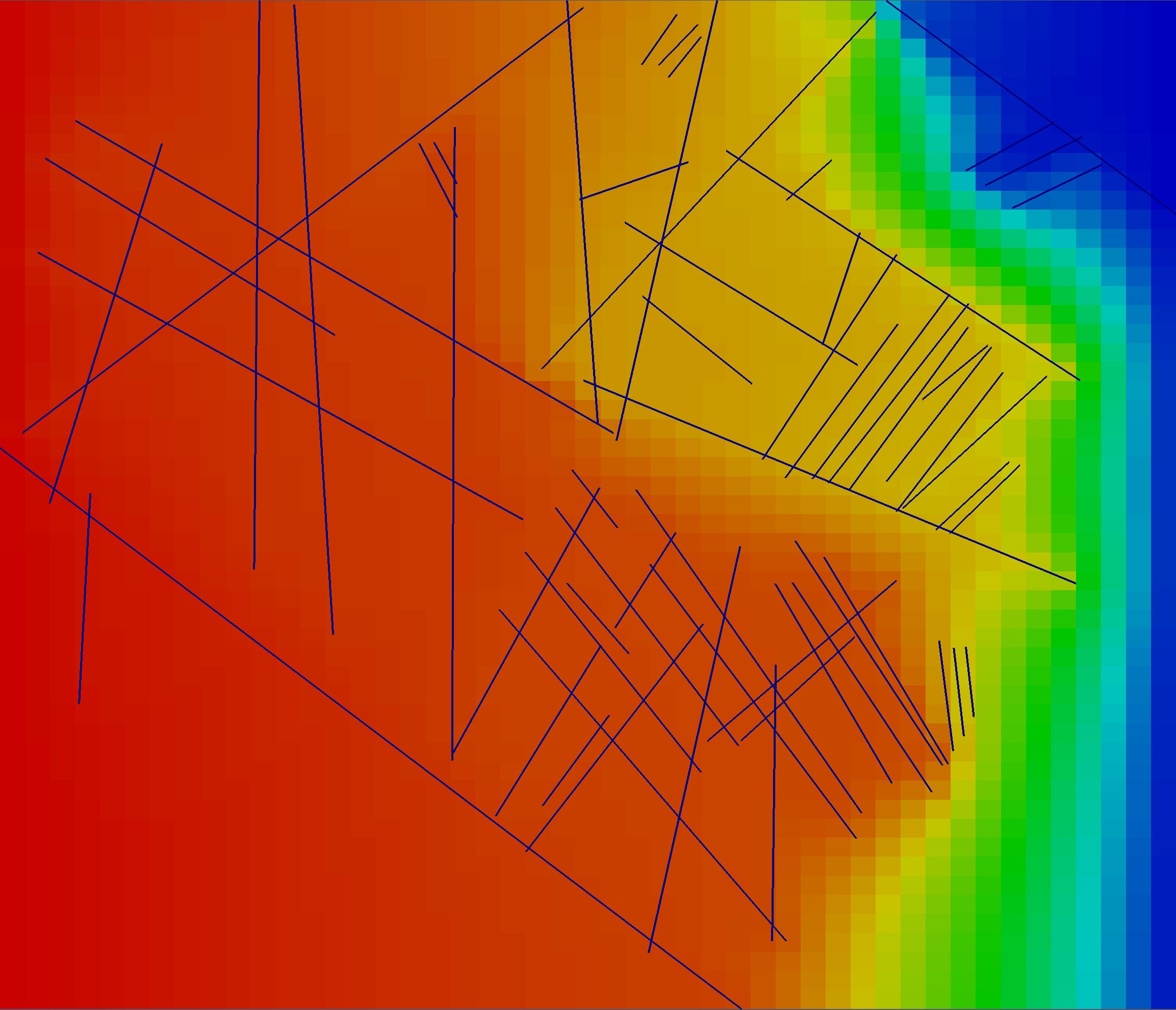}%
    }%
    \subfloat[mortar-DFM]{
    \includegraphics[width=0.5\textwidth]{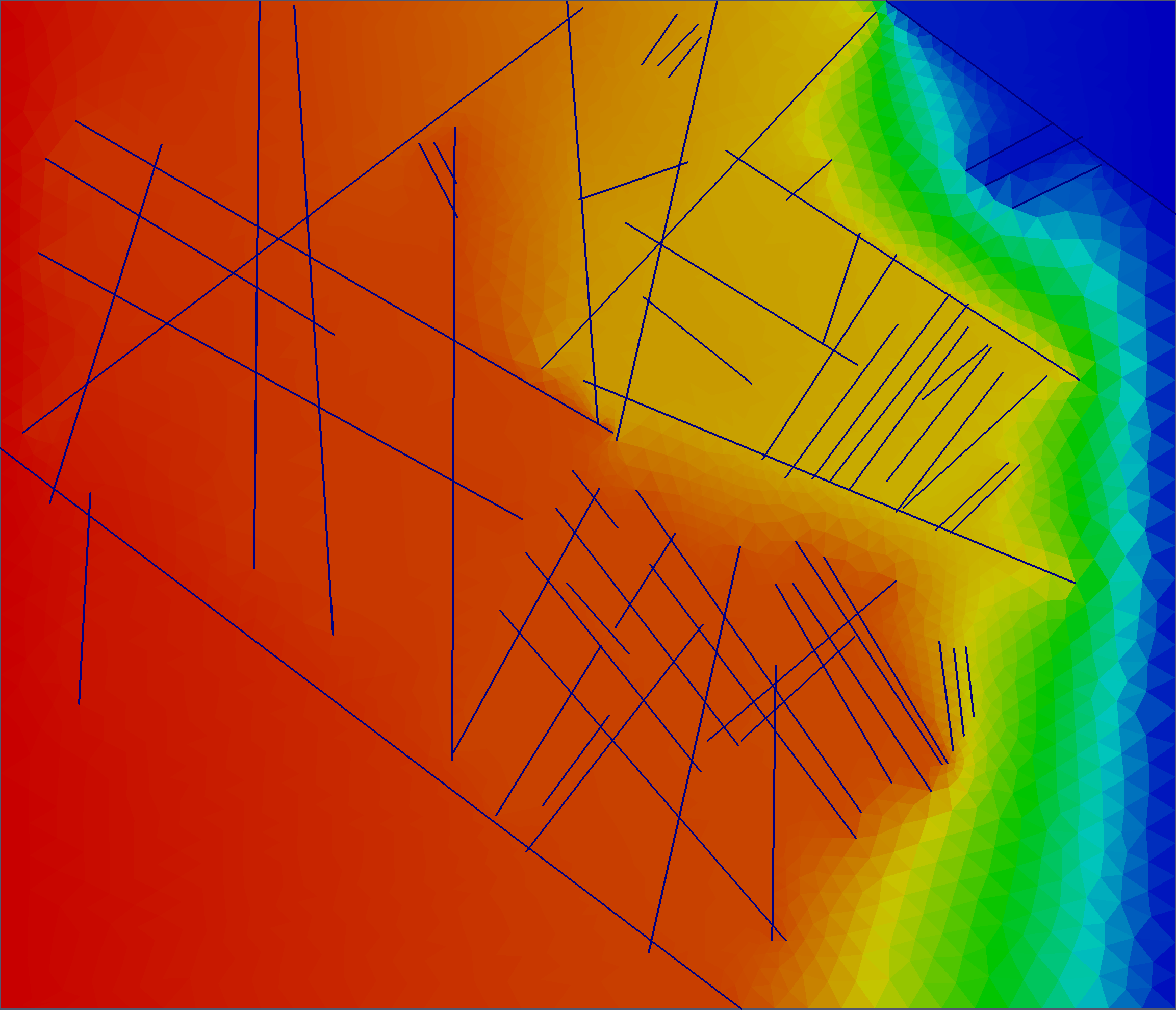}%
    }%
    \caption{Representation of the matrix pressures field for the realistic case.
    The solution values range between $0$ and $\unit[101325]{Pa}$.}%
    \label{fig:solution_complex}
\end{figure}

\begin{figure}
    \centering
    \subfloat[Box-DFM]{
    \includegraphics[width=0.5\textwidth]{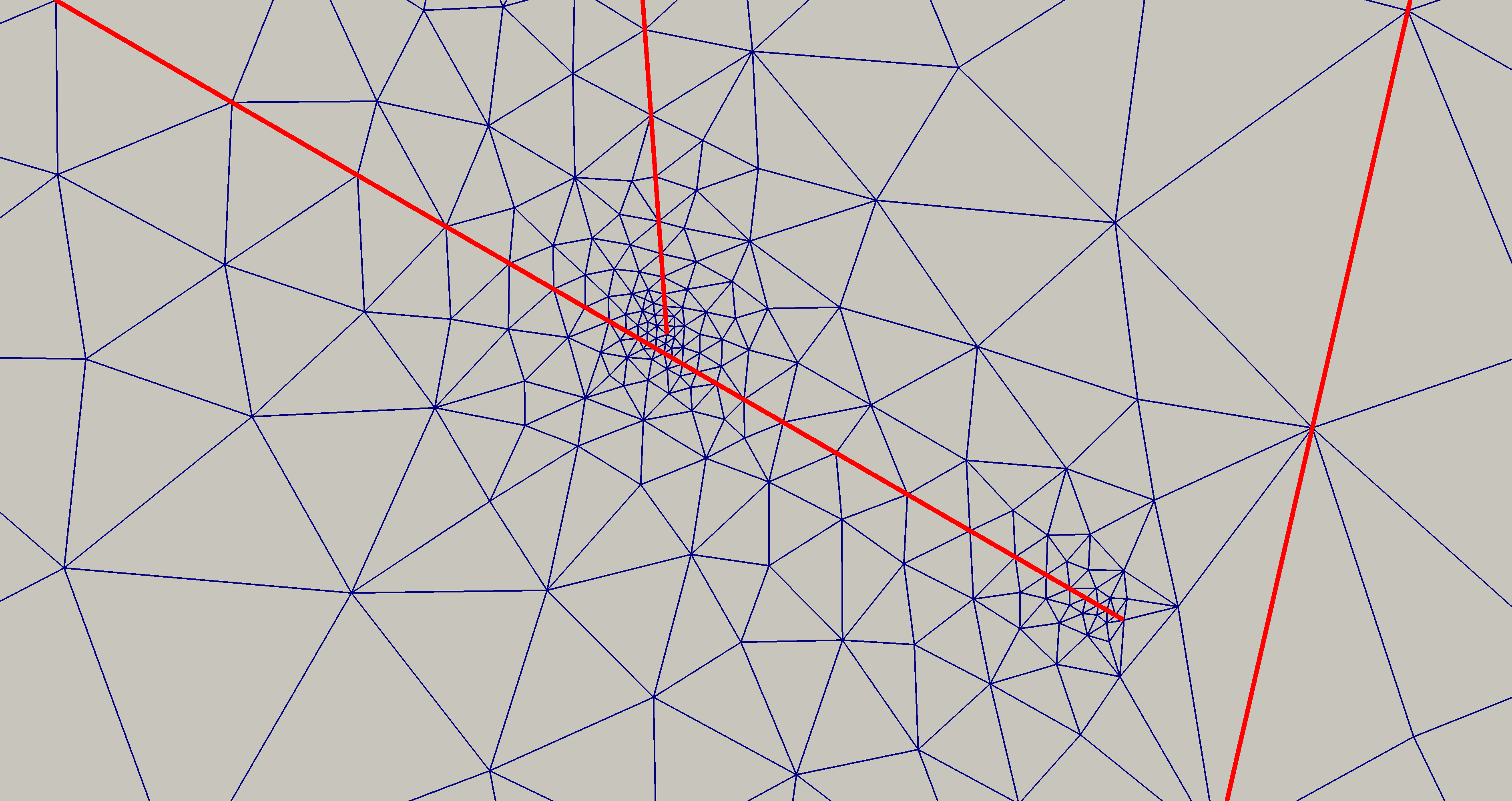}%
    }
    \subfloat[CC-DFM]{
    \includegraphics[width=0.5\textwidth]{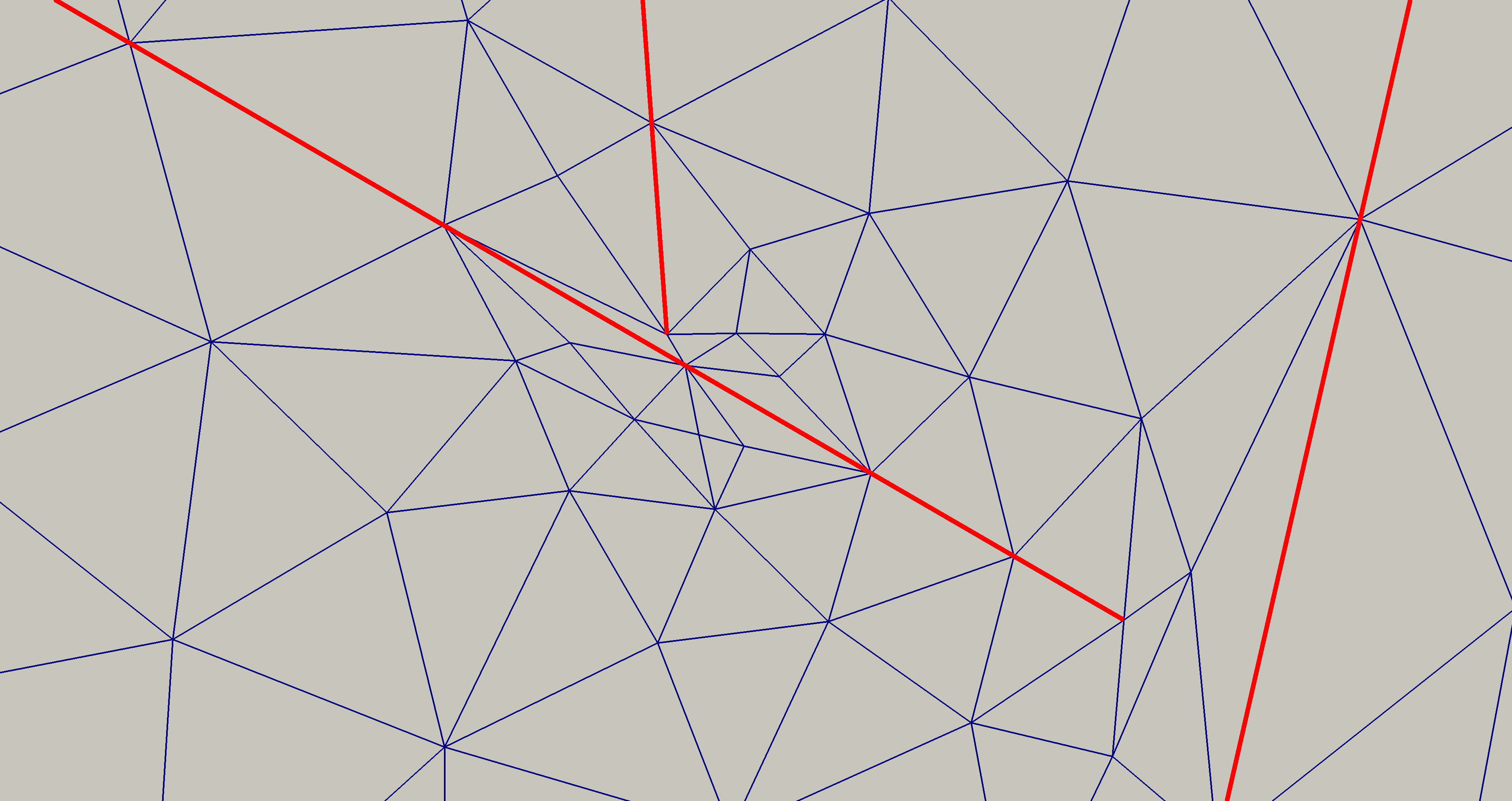}
    }\\
    \subfloat[EDFM]{
    \includegraphics[width=0.5\textwidth]{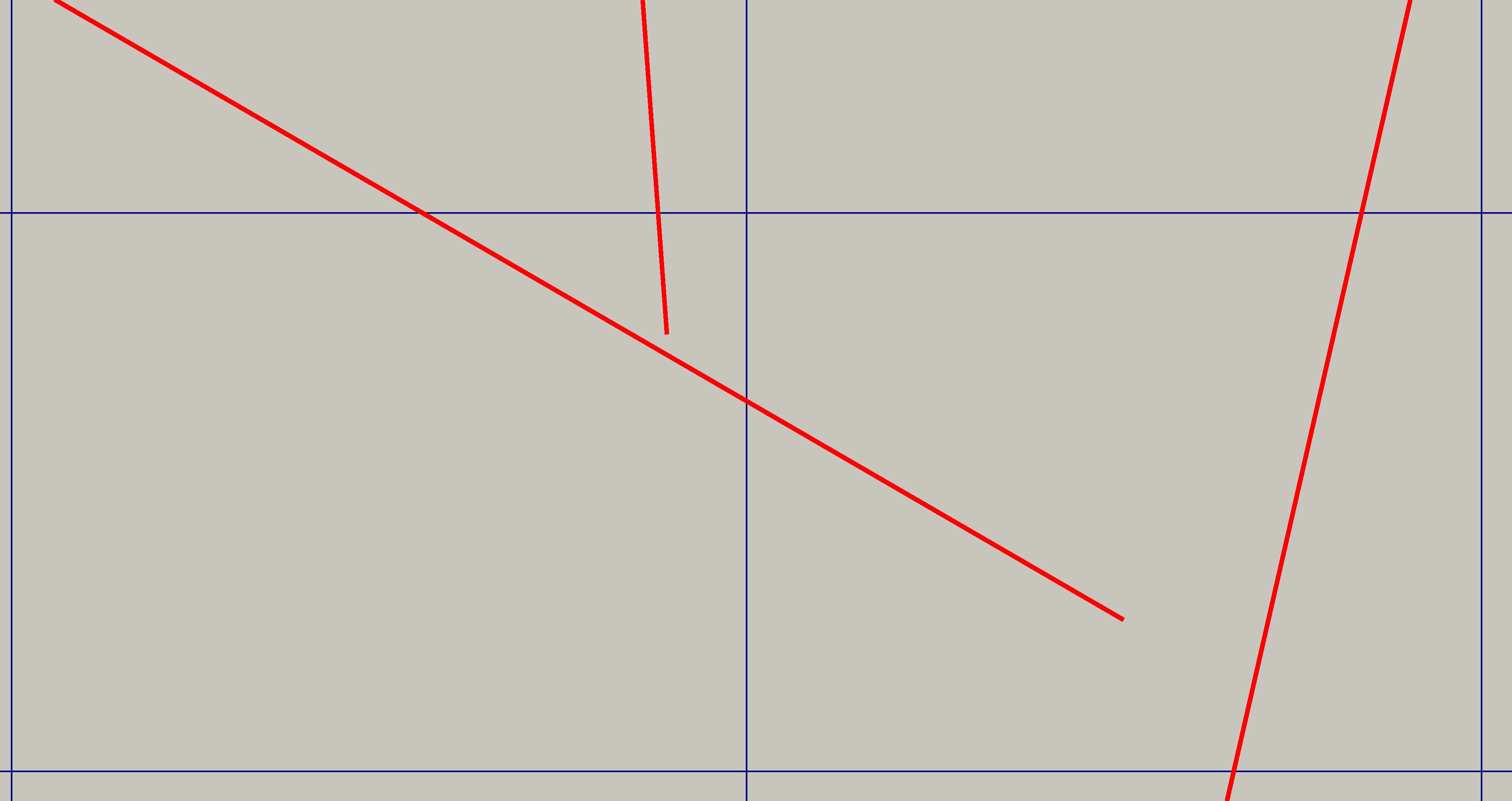}%
    }%
    \subfloat[mortar-DFM]{
    \includegraphics[width=0.5\textwidth]{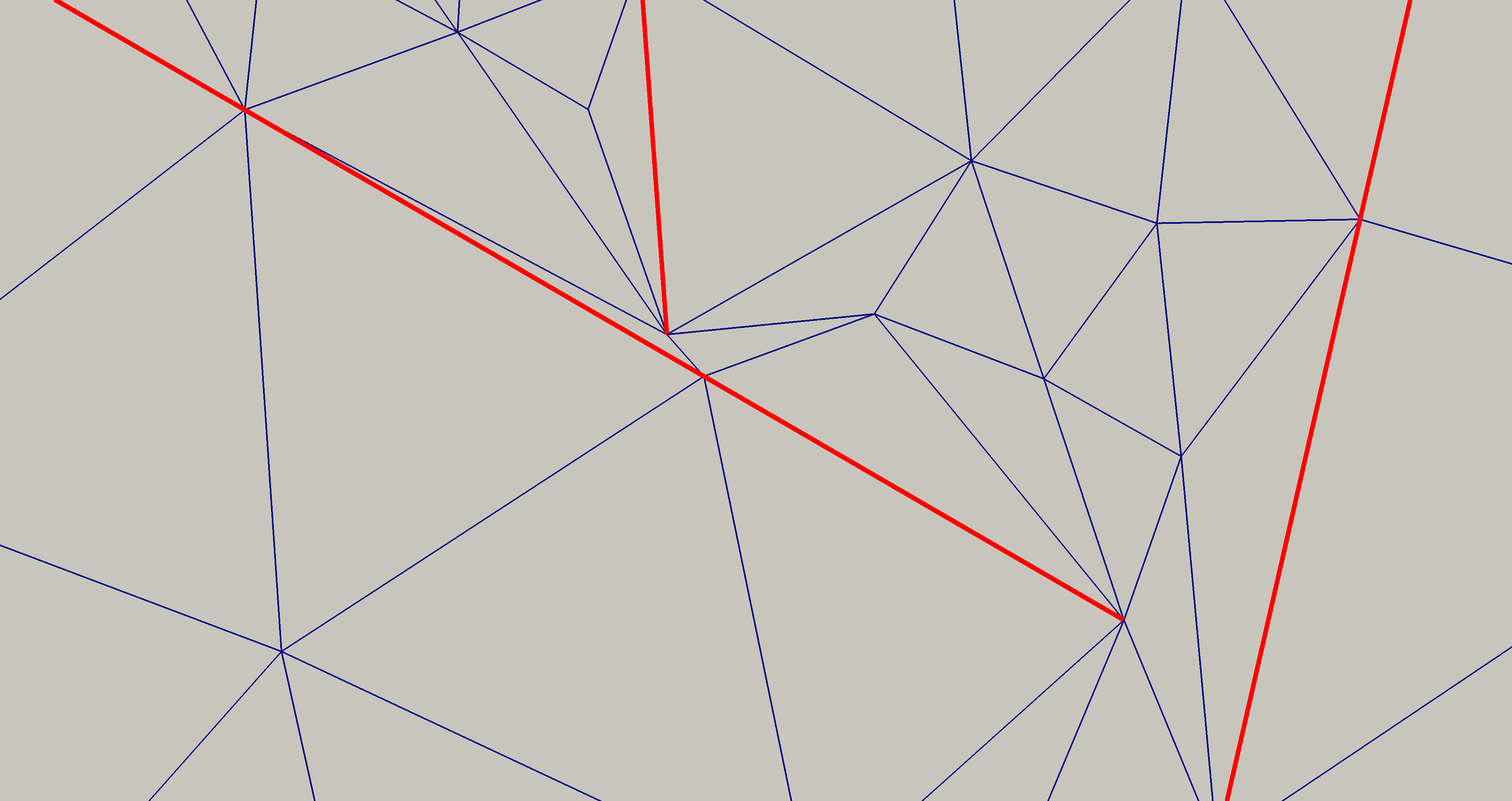}%
    }%
    \caption{Benchmark 4: Representation of mesh in the middle of the domain. The size of the
    picture is approximately $\unit[30]{m} \times \unit[15]{m}$ centered in $(360,350)$.
    It is represented by the small rectangle in the centre of Figure
    \ref{fig:sotra} left.}%
    \label{fig:solution_complex_zoom}
\end{figure}

\begin{figure}
    \centering
    \subfloat[$y=500\,\text{m}$]{
    \includegraphics[width=0.5\textwidth]{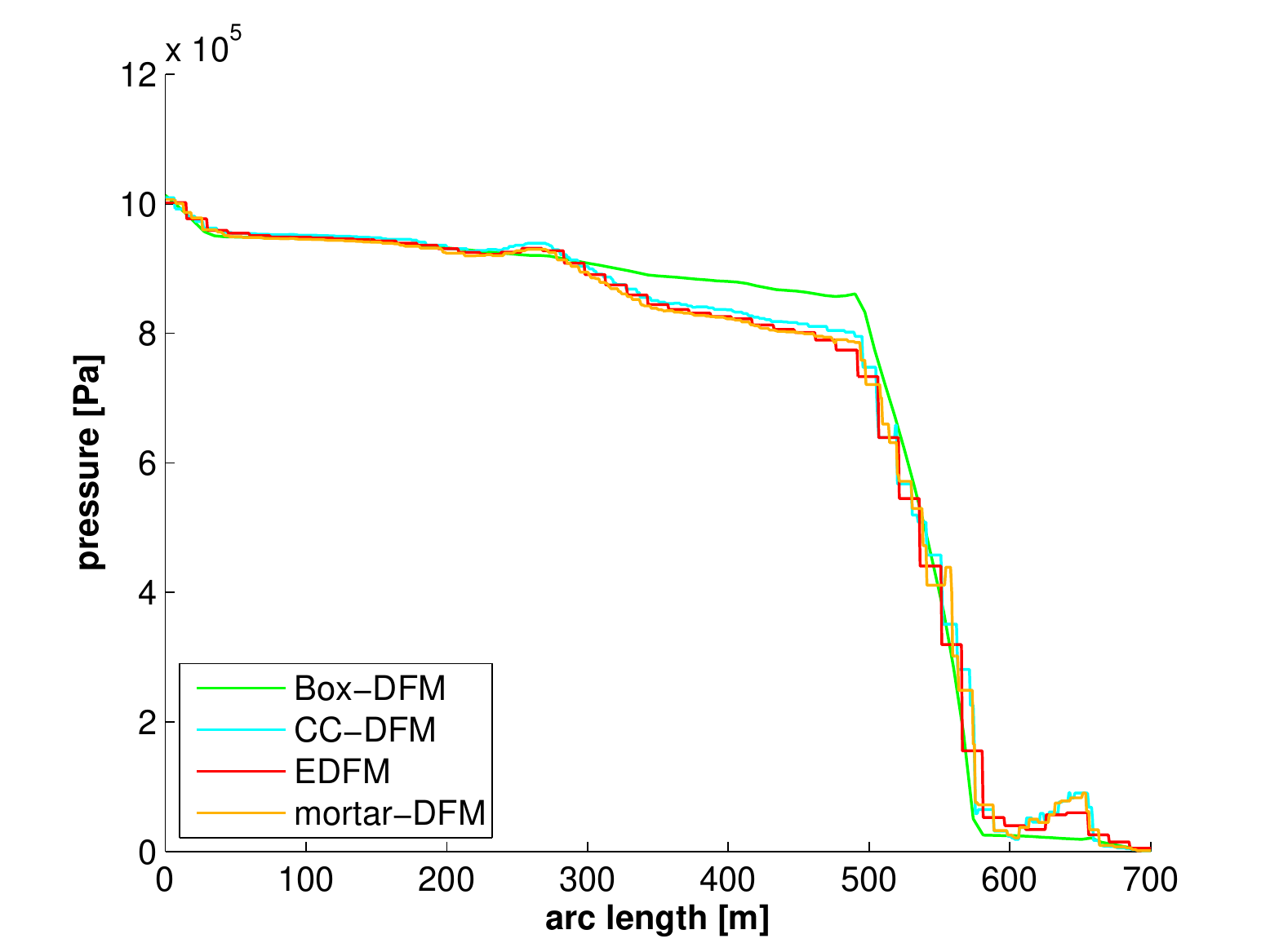}%
    }
    \subfloat[$x=625\,\text{m}$]{
    \includegraphics[width=0.5\textwidth]{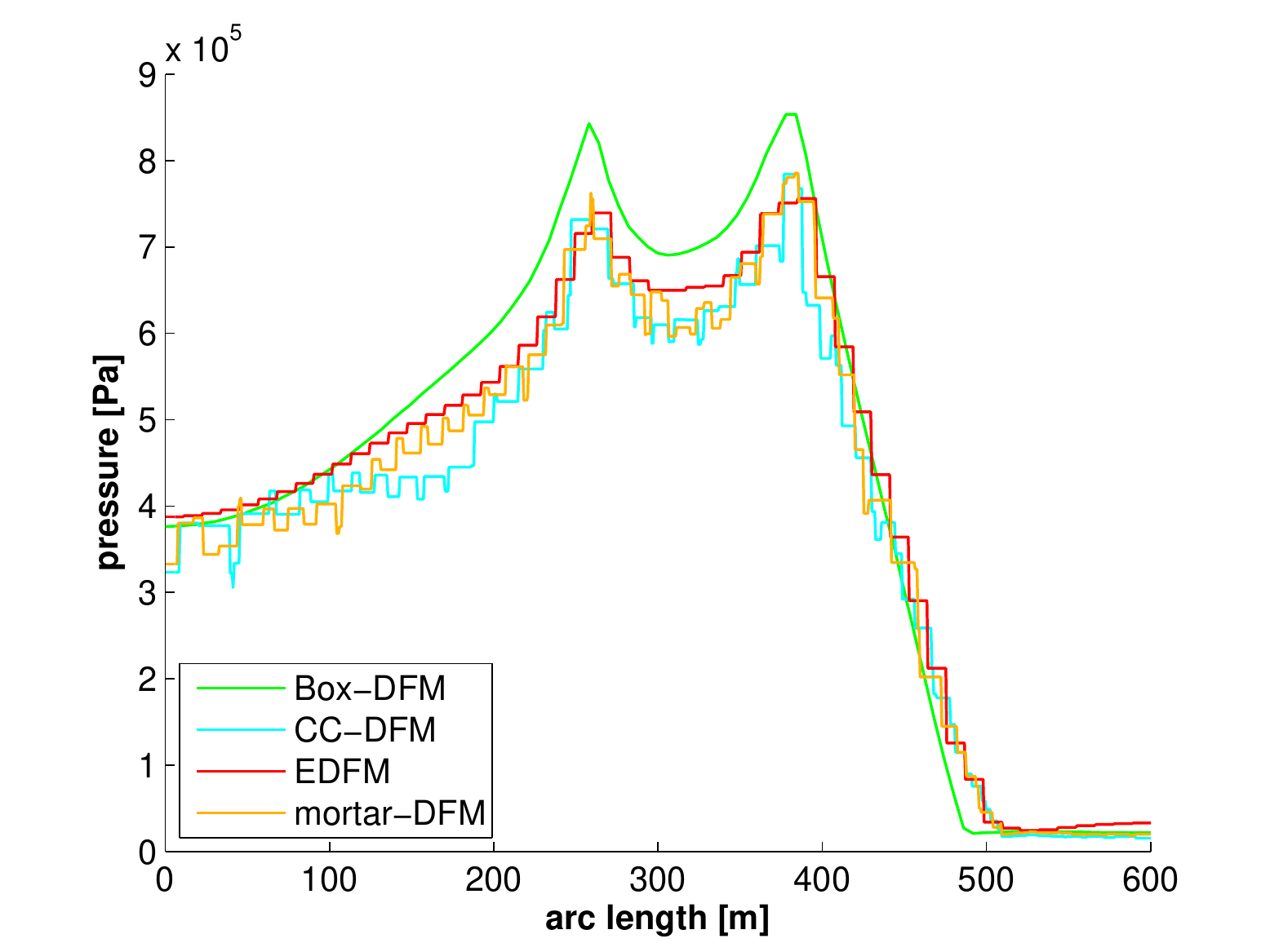}%
    }%
    \caption{Benchmark 4: Pressure solutions of the 4 methods plotted over lines
    (a) $y=\unit[500]{m}$, and (b) $x = \unit[625]{m}$.}%
    \label{fig:solution_complex_line}
\end{figure}

Table \ref{tbl:complex_error} lists the number of degrees of freedom, the
density of the associated matrix, and its condition number for the different
methods.
\renewcommand{\arraystretch}{1.1}
\begin{table}[hbt]
\centering
\begin{tabular}{|l|c|c|c|c|c|}\hline
\textbf{method} & \textbf{d.o.f.} & \textbf{matrix elem} &
\textbf{frac elem}
& \textbf{nnz/size$^2$} & $\|\cdot\|_2$\textbf{-cond} \\\hline
Box-DFM    & 5563 & 10807 triangles & 1386 & 1.2e-3 & 9.3e5 \\\hline
CC-DFM     & 8481 & 7614 triangles & 867 & 4.9e-4 & 5.3e6\\\hline
EDFM       & 3599  & 2491 quads & 1108 & 1.4e-3  & 4.7e6 \\\hline
mortar-DFM & 25258 & 8319 triangles & 1317 & 2.0e-4 & 2.2e17\\\hline
\end{tabular}
\caption{Discretization and matrix characteristics for Benchmark \theBenchmarkCounter.}
\label{tbl:complex_error}
\end{table}
\renewcommand{\arraystretch}{1.0}
Due to the geometrical difficulties of the network the request of having a
similar number of degrees of freedom among the methods is relaxed, as Table
\ref{tbl:complex_error} indicates. Considering Figure
\ref{fig:solution_complex}, the solutions are reported for the four methods. We
notice that, except for the top right part of the domain in the Box-DFM method,
the solutions are similar and comparable which is an indication of their
correctness.
Compared to the previous test cases the mesh generation is the main concern and
some of the methods require a fine tuning to avoid unphysical connections among
elements where the fracture are close. An example can be found in the middle of
the domain and reported in Figure \ref{fig:solution_complex_zoom}. Only EDFM is
more robust with respect to this constraint.
To present a more detailed comparison among the methods, Figure
\ref{fig:solution_complex_line} represents the pressure solution along two different
lines: for $y=\unit[500]{m}$ and for $x=\unit[625]{m}$. We note that the methods behave
similarly, and the Box-DFM slightly overestimates some peaks. The oscillation
of the methods are related to mesh effects.

\clearpage
\section{Summary and Outlook}
\label{sec:summary}
Four benchmark cases for single-phase flow in fractured porous media have been proposed
and employed to compare the performances of several state-of-the-art hybrid-dimensional discrete-fracture-matrix
models.
If we consider the cases where all the methods are employed within the applicability range
for which they were originally developed,
the results are in quite good agreement. In particular, fracture networks exhibiting a larger permeability
than the surrounding matrix can be accurately described by all methods.
On the other hand, not all methods are capable of modeling blocking fractures.
In this case some methods fail to predict the correct flow patterns for the corresponding cases.
Especially noteworthy are the large differences in the condition numbers of the associated
system matrices. The effect of these differences on the behavior of linear solvers is difficult
to quantify in a comparable manner, since the different methods pose different requirements
for such solvers.
In principle, all participating methods should have been able to run all proposed cases.
However, due to implementation restrictions, some methods could not perform the cases with
more complex fracture network geometries.
Even if the methodology is general enough, technical difficulties can become
crucial obstacles to tackling realistic scenarios.

All the investigated benchmarks are restricted to simple physics and two-dimensional computational
domains. This should give other researchers developing DFM models the chance to perform
comparison studies for their methods. We encourage the scientific community to contribute their
results for the benchmarks to a corresponding Git repository at
\url{https://git.iws.uni-stuttgart.de/benchmarks/fracture-flow}.

Further benchmark cases may be developed in the near future.
In particular, we are very interested in enhancing the
purely single-phase single-component flow physics by adding
transport, deformation and/or reaction processes.
We aim to carry out these efforts in a broader context by means of international workshops.

 \appendix
\section{Acknowledgement}

The authors warmly thank Luisa F. Zuluaga, from University of Bergen,
for constructing and providing the real fracture network for the example in
Section \ref{sec:realistic}. The authors wish to thank also Luca Pasquale and
Stefano Zonca.

The second author acknowledges financial support from the GeoStim project from the Research Council of Norway (project no. 228832) through the ENERGIX program.
The third author was supported by Norwegian Research Council grant 233736.
The fourth author acknowledges financial support from the ANIGMA project
from the Research Council of Norway (project no. 244129/E20) through the ENERGIX program.

\section{Domain modifications for Benchmark 1}
\label{sec:hydrocoin_mod}
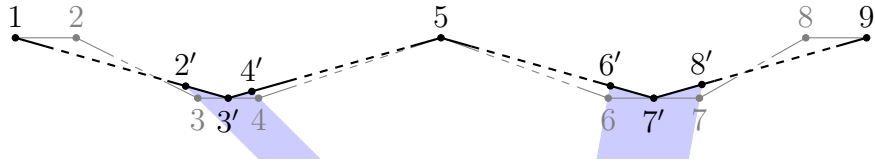
\begin{figure}[hbt]
\begin{center}
\begin {tikzpicture}[scale=0.8]

\draw[blue!20, fill=blue!20] (2.8, 0.2) -- (3.5, 0) -- (3.888888888888889, 0.111111111111111) -- (5, -1) -- (4, -1) -- cycle;
\draw[blue!20, fill=blue!20] (9.787037037037036, 0.203703703703704) -- (10.5, 0) -- (11.291095890410960, 0.226027397260274) -- (11.0681818182, -1) -- (9.5681818182, -1) -- cycle;

\draw[gray] (0,1) -- (1, 1) -- (1.5, 0.75);
\draw[dashed, gray] (1.5, 0.75) -- (2.5, 0.25);
\draw[gray] (2.5, 0.25) -- (3, 0) -- (4, 0) -- (4.5, 0.166667);
\draw[dashed, gray] (4.5, 0.166667) -- (6.5, 0.833333);
\draw[gray] (6.5, 0.833333) -- (7, 1) -- (7.5, 0.81818181818);
\draw[dashed, gray] (7.5, 0.81818181818) -- (9.25, 0.18181818);
\draw[gray] (9.25, 0.18181818) -- (9.75, 0) -- (11.25, 0) -- (11.75, 0.28571428571);
\draw[dashed, gray] (11.75, 0.28571428571) -- (12.5, 0.71428571428);
\draw[gray] (12.5, 0.71428571428) -- (13, 1) -- (14, 1);

\draw[thick] (0,1) -- (0.5, 0.85714285714);
\draw[thick, dashed] (0.5, 0.85714285714) -- (2.5, 0.28571428571);
\draw[thick] (2.5, 0.28571428571) -- (3.5, 0) -- (4.5, 0.28571428571);
\draw[thick, dashed] (4.5, 0.28571428571) -- (6.5, 0.85714285714);
\draw[thick] (6.5, 0.85714285714) -- (7, 1) -- (7.5, 0.85714285714);
\draw[thick, dashed] (7.5, 0.85714285714) -- (9.5, 0.28571428571);
\draw[thick] (9.5, 0.28571428571) -- (10.5, 0) -- (11.5, 0.28571428571);
\draw[thick, dashed] (11.5, 0.28571428571) -- (13.5, 0.85714285714);
\draw[thick] (13.5, 0.85714285714) -- (14, 1);

\filldraw[black] (0, 1) circle (1.5pt) node[above] {1};
\filldraw[gray] (1, 1) circle (1.5pt) node[above] {2};
\filldraw[gray] (3, 0) circle (1.5pt) node[below] {3};
\filldraw[gray] (4, 0) circle (1.5pt) node[below] {4};
\filldraw[black] (7, 1) circle (1.5pt) node[above] {5};
\filldraw[gray] (9.75, 0) circle (1.5pt) node[below] {6};
\filldraw[gray] (11.25, 0) circle (1.5pt) node[below] {7};
\filldraw[gray] (13, 1) circle (1.5pt) node[above] {8};
\filldraw[black] (14, 1) circle (1.5pt) node[above] {9};

\filldraw[black] (2.8, 0.2) circle (1.5pt) node[above] {2$^\prime$};
\filldraw[black] (3.5, 0) circle (1.5pt) node[below] {3$^\prime$};
\filldraw[black] (3.888888888888889, 0.111111111111111) circle (1.5pt) node[above] {4$^\prime$};
\filldraw[black] (9.787037037037036, 0.203703703703704) circle (1.5pt) node[above] {6$^\prime$};
\filldraw[black] (10.5, 0) circle (1.5pt) node[below] {7$^\prime$};
\filldraw[black] (11.291095890410960, 0.226027397260274) circle (1.5pt) node[above] {8$^\prime$};

\end {tikzpicture}
\end{center}
\caption{Modifications of the Hydrocoin model domain compared to the original
formulation \cite{hydrocoin}. The original upper boundary is drawn with gray
thin lines, while thick black lines are used for the modified boundary. Modified
node locations are indicated by numbers superscripted with $\mbox{}^\prime$. The
shaded regions show the upper parts of the two slightly extended
equi-dimensional fractures.}
\label{fig:hydrocoin_modifications}
\end{figure}

Table \ref{tbl:hydrocoin_setup} provides the exact coordinates of 
the points from Figure \ref{fig:hydrocoin_setup}.
\begin{table}[hbt]
\centering
\caption{Coordinates of the numbered points in the modeled region of the problem
depicted in Figure \ref{fig:hydrocoin_setup}.}
\begin{tabular}{rrr|rrr}
\toprule
pt & $x$ (m) & $z$ (m) & pt & $x$ (m) & $z$ (m) \\
\midrule
 1 &    0 & 150 & 11 & 1505  & -1000  \\
 2$^\prime$ &   394.285714286 & 100.714285714 & 12 & 1495  & -1000  \\
 3$^\prime$ &  400 & 100 & 13  & 1007.5  & -1000  \\
 4$^\prime$ &  404.444444444 & 100.555555556 & 14 &  992.5 &  -1000 \\
 5 &  800 & 150 & 15 &  0 &  -1000 \\
 6$^\prime$ & 1192.66666667 & 100.916666667 & 16 & 1071.34615385  &  -566.346153846 \\
 7$^\prime$ & 1200 & 100 & 17 & 1084.03846154  &  -579.038461538 \\
 8$^\prime$ & 1207.6744186 & 100.959302326 & 18 & 1082.5  & -587.5  \\
 9 & 1600 & 150 & 19 & 1069.80769231  &  -574.807692308 \\
10 & 1600 & -1000 &  &   &   \\
\bottomrule
\end{tabular}
\label{tbl:hydrocoin_setup}
\end{table}
In comparison to the original setup, the plateaus close to the upper
left and right corners 1 and 9 have been omitted. Moreover, the upper ends of
the two fractures have been modified according to
Figure \ref{fig:hydrocoin_modifications} which amounts to the changes of nodes
2--4 and 6--8.
Finally, the position of nodes 16--19 has been recalculated with higher
precision. The hybrid-dimensional models do not take into account nodes 2,4,6,8
and 16--19 and combine nodes 11,12 and 13,14, since the two-dimensional fracture
regions have been reduced to two intersecting straight lines.

\section{Fracture coordinates for Benchmark 3}
\label{sec:anna_coord}
The coordinates are listed in Table \ref{tbl:complex_domain_anna}.
\begin{table}[hbt]
\centering
 \begin{tabular}{|c|c|c|c|c|}
 \hline
      Nf  &      xA    &    yA   &     xB    &    yB\\
\hline
    1   &      0.0500  &  0.4160  &  0.2200  &  0.0624\\
    2   &      0.0500  &  0.2750  &  0.2500  &  0.1350\\
    3   &      0.1500  &  0.6300  &  0.4500  &  0.0900\\
    4   &      0.1500  &  0.9167  &  0.4000  &  0.5000\\
    5   &      0.6500  &  0.8333  &  0.8500  &  0.1667\\
    6   &      0.7000  &  0.2350  &  0.8500  &  0.1675\\
    7   &      0.6000  &  0.3800  &  0.8500  &  0.2675\\
    8   &      0.3500  &  0.9714  &  0.8000  &  0.7143\\
    9   &      0.7500  &  0.9574  &  0.9500  &  0.8155\\
   10   &      0.1500  &  0.8363  &  0.4000  &  0.9727\\
   \hline
 \end{tabular}
\caption{Benchmark 3: Fracture coordinates}%
\label{tbl:complex_domain_anna}
\end{table}

\ifpreprint
\else
\section*{Bibliography}
\fi
\bibliographystyle{elsarticle-harv}
\bibliography{literature}

\end{document}